%% file: User-Friendly-focm-vf.tex
\documentclass[11pt,letterpaper,twoside,reqno,draft]{amsart}

\input{macro-file}

\usepackage{subfigure}
\usepackage{wasysym}

\begin{document}
\title[Tail Bounds for Sums of Random Matrices]
{User-Friendly Tail Bounds \\ for Sums of Random Matrices}

\author{Joel A.~Tropp}

\keywords{Discrete-time martingale, large deviation, probability inequality, random matrix, sum of independent random variables}

\thanks{2010 {\it Mathematics Subject Classification}.
Primary:
60B20. 
Secondary:
60F10, 
60G50, 
60G42
}

\thanks{JAT is with Computing \& Mathematical Sciences, MC 305-16, California Institute of Technology, Pasadena, CA 91125.
E-mail: \url{jtropp@cms.caltech.edu}.
Research supported by ONR award N00014-08-1-0883, DARPA award N66001-08-1-2065, and AFOSR award FA9550-09-1-0643.}

\date{25 April 2010.  Revised on 15 June 2010, 14 November 2010, 7 January 2011, and 15 June 2011.}

\begin{abstract}
This paper presents new probability inequalities for sums of independent, random, self-adjoint matrices.  These results place simple and easily verifiable hypotheses on the summands, and they deliver strong conclusions about the large-deviation behavior of the maximum eigenvalue of the sum.  Tail bounds for the norm of a sum of random rectangular matrices follow as an immediate corollary.  The proof techniques also yield some information about matrix-valued martingales.

In other words, this paper provides noncommutative generalizations of the classical bounds associated with the names Azuma, Bennett, Bernstein, Chernoff, Hoeffding, and McDiarmid.  The matrix inequalities promise the same diversity of application, ease of use, and strength of conclusion that have made the scalar inequalities so valuable.
\end{abstract}



\maketitle

\section{Introduction}



Random matrices have come to play a significant role in computational mathematics.  This line of research has advanced by using established methods from random matrix theory, but it has also generated difficult questions that cannot be addressed without new tools.  Let us summarize some of the challenges that arise in numerical applications.



\begin{itemize} \setlength{\itemsep}{0.5pc}
\item	Research has extended well {beyond the classical ensembles} (e.g., Wishart matrices and Wigner matrices) to encompass many other classes of random matrices.  For instance, it is now common to study the properties of a sparse matrix sampled from a fixed matrix or a random submatrix drawn from a fixed matrix.

\item	We also encounter highly structured matrices that involve a {limited amount of randomness}.  One important example is the randomized DFT, which consists of a diagonal matrix of random signs multiplied by a discrete Fourier transform matrix. 

\item	Questions about the spectral properties of random matrices remain fundamental, but modern problems can also involve other considerations.  For example, we might need to estimate the cut norm of a random adjacency matrix.  Or we might want to study the action of a random operator on a class of vectors or matrices.

\item	Most problems in numerical mathematics concern matrices of {finite order}.  Asymptotic theory is less relevant in practice.  

\item	We often require explicit {large-deviation theorems} for statistics of random matrices so that we can study rates of convergence. 

\item	Results with {effective constants} are essential to ensure that algorithms are provably correct.
\end{itemize}

We have encountered these issues in a wide range of problems from computational mathematics:
smoothed analysis of Gaussian elimination~\cite{SST06:Smoothed-Analysis};
semidefinite relaxation and rounding of quadratic maximization problems~\cite{Nem07:Sums-Random,So09:Moment-Inequalities};
construction of maps for dimensionality reduction~\cite{AC09:Fast-Johnson-Lindenstrauss};
matrix approximation by sparsification~\cite{AM07:Fast-Computation} and by sampling submatrices~\cite{RV07:Sampling-Large};
analysis of sparse approximation~\cite{Tro08:Conditioning-Random}
and compressive sampling~\cite{CR07:Sparsity-Incoherence} algorithms;
randomized schemes for low-rank matrix factorization~\cite{HMT11:Finding-Structure};
and analysis of algorithms for completion of low-rank matrices~\cite{Gro11:Recovering-Low-Rank,Rec09:Simpler-Approach}.
%
And this list is by no means comprehensive!

In most of these applications, the methods currently invoked to study random matrices require a substantial amount of practice to use effectively.  Even so, the final results tend to be a little disappointing: the constants are usually poor and the predictions are sometimes coarser than we might like.  These frustrations have led us to search for simpler techniques that still yield detailed quantitative information about finite random matrices.

\subsection{Technical Overview}

We consider a finite sequence $\{\mtx{X}_k\}$ of 
random, self-adjoint matrices with dimension $d$.  Our goal is to harness basic properties of these matrices to bound the probability
\begin{equation}\label{eqn:intro-prob}
\Prob{ \lambda_{\max}\left( \sum\nolimits_k \mtx{X}_k \right) \geq t }.
\end{equation}
Here and elsewhere, $\lambda_{\max}$ denotes the algebraically largest eigenvalue of a self-adjoint matrix.  This formulation is more general than it may appear because we can exploit the same ideas to explore several related problems:
\begin{itemize}
\item  We can study the smallest eigenvalue of the sum.
\item	We can bound the largest singular value of a sum of random rectangular matrices.
\item	Related arguments apply to matrix martingales and other adapted sequences.
\end{itemize}
Indeed, the expression~\eqref{eqn:intro-prob} captures the essence of many questions that arise in numerical applications of random matrix theory, including most of the research cited above.

Observe that~\eqref{eqn:intro-prob} formally resembles the probability that a sum of real random variables exceeds a certain level.  
The \term{Laplace transform method}, attributed to Bernstein, is a particularly elegant system for producing tail bounds for sums of scalar random variables; see~\cite{McD98:Concentration,Lug09:Concentration-Measure} for accessible discussions.  In a remarkable paper~\cite{AW02:Strong-Converse}, Ahlswede and Winter show how to transport the Laplace transform method to the matrix setting.  They establish that
\begin{equation} \label{eqn:intro-lt}
\Prob{ \lambda_{\max}\left( \sum\nolimits_k \mtx{X}_k \right) \geq t }
	\leq \inf_{\theta > 0} \left\{ \econst^{-\theta t} \cdot
	\Expect \trace \exp\left( \sum\nolimits_k \theta \mtx{X}_k \right) \right\}.
\end{equation}
In words, the probability~\eqref{eqn:intro-prob} is controlled by a matrix version of the moment generating function (mgf).  See Proposition~\ref{prop:laplace-transform} for an easy proof of~\eqref{eqn:intro-lt} that is due to Oliveira~\cite{Oli10:Sums-Random}.

The matrix Laplace transform estimate~\eqref{eqn:intro-lt} presents a serious technical challenge.  We must control the trace of the matrix mgf
$$
\Expect \trace \exp\left( \sum\nolimits_k \theta \mtx{X}_k \right)
$$
using information about the summands $\mtx{X}_1, \mtx{X}_2, \mtx{X}_3, \dots$.  This estimate requires powerful tools, and it stands as the major impediment to bounding the tail probability~\eqref{eqn:intro-prob}. 

The true significance of the Ahlswede--Winter argument~\cite[App.]{AW02:Strong-Converse} consists in their technique for computing the required bounds on the matrix mgf.  We describe their method in~\S\ref{sec:aw}.  The following probability inequality for a matrix Gaussian series is typical of the results that emerge from their approach.  Let $\{ \mtx{A}_k \}$ be a family of fixed self-adjoint matrices with dimension $d$, and let $\{\gamma_k\}$ be a sequence of independent standard normal variables.  Then
\begin{equation} \label{eqn:aw-intro}
\Prob{ \lambda_{\max}\left( \sum\nolimits_k \gamma_k \mtx{A}_k \right) \geq t }
	\leq d \cdot \econst^{-t^2/2\sigma_{\rm AW}^2}
	\quad\text{where}\quad
	\sigma_{\rm AW}^2 := \sum\nolimits_k \lambda_{\max}( \mtx{A}_k^2 ).
\end{equation}
The Ahlswede--Winter apparatus leads to a collection of other interesting probability inequalities; see~\S\ref{sec:related} for references.  Nevertheless, tail bounds developed in this fashion, including~\eqref{eqn:aw-intro}, are usually very far from optimal.  See~\S\ref{sec:aw} and~\S\ref{sec:aw-gauss} for further discussion of this point.




This paper describes a more satisfactory framework for completing the bound on the matrix mgf.
The crucial new ingredient in our argument is a deep theorem~\cite[Thm.~6]{Lie73:Convex-Trace} of Lieb from his seminal paper on convex trace functions. 
We introduce Lieb's theorem in~\S\ref{sec:lieb}, and we explain how to combine this result with the matrix Laplace transform technique.  We use this scheme to obtain a large family of probability inequalities that are essentially sharp in a wide variety of situations.  

Our approach represents a dramatic advance beyond the Ahlswede--Winter technique.  For example, our method delivers the following bound for a matrix Gaussian series:
\begin{equation} \label{eqn:gauss-intro}
\Prob{ \lambda_{\max}\left( \sum\nolimits_k \gamma_k \mtx{A}_k \right) \geq t }
	\leq d \cdot \econst^{-t^2/2\sigma^2}
	\quad\text{where}\quad
	\sigma^2 := \lambda_{\max}\left( \sum\nolimits_k  \mtx{A}_k^2 \right).
\end{equation}
The estimate~\eqref{eqn:gauss-intro} offers a fundamental advantage over~\eqref{eqn:aw-intro} because the variance parameter $\sigma^2$ is often $d$ times smaller than $\sigma^2_{\rm AW}$.  Furthermore, the discussion in~\S\ref{sec:rademacher} demonstrates that the inequality~\eqref{eqn:gauss-intro} cannot be sharpened without changing its structure.  This improvement is typical of results constructed from our blueprint.

\subsection{Index of Inequalities}

This work contains a large number of bounds for the probability~\eqref{eqn:intro-prob}.  The precise form of each inequality depends on prior information about the summands.  As a service to the reader, we have collected the most useful results in this section.  We have also included a short qualitative discussion of each bound, along with the location in the paper where the  full treatment appears.


\subsubsection{Notation}

The symbol $\psdle$ denotes the semidefinite order on self-adjoint matrices.  The maps $\lambda_{\min}$ and $\lambda_{\max}$ return the algebraically smallest and largest eigenvalue of a self-adjoint matrix.  We write $\norm{ \cdot }$ for the spectral norm, which equals the largest singular value of a matrix.

\subsubsection{Main Results for Positive-Semidefinite Matrices}

In classical probability theory, one of the most famous concentration results concerns the number of successes in a sequence of independent random trials.  This quantity can be expressed as a sum of independent, bounded random variables.  Chernoff's large-deviation theorem~\cite{Che52:Measure-Asymptotic} provides explicit estimates on the probability that this type of series is greater than (or smaller than) a specified level.


In the matrix setting, the analogous theorem concerns a sum of positive-semidefinite random matrices subject to a uniform eigenvalue bound.  The matrix Chernoff inequality shows that the extreme eigenvalues of the matrix series have the same binomial-type behavior that occurs in the scalar case.


\begin{thm}[Matrix Chernoff]
Consider a finite sequence $\{ \mtx{X}_k \}$ of independent, random, self-adjoint matrices with dimension $d$.  Assume that each random matrix satisfies
$$
\mtx{X}_k \psdge \mtx{0}
\quad\text{and}\quad
\lambda_{\max}( \mtx{X}_k ) \leq R
\quad\text{almost surely}.
$$
Define
$$
\mu_{\min} := \lambda_{\min}\left( \sum\nolimits_k \Expect \mtx{X}_k \right)
\quad\text{and}\quad
\mu_{\max} := \lambda_{\max}\left( \sum\nolimits_k \Expect \mtx{X}_k \right).
$$
Then
\begin{align*}
\Prob{ \lambda_{\min}\left( \sum\nolimits_k \mtx{X}_k \right) \leq (1 - \delta) \mu_{\min} }
	&\leq d \cdot \left[ \frac{\econst^{-\delta}}{(1-\delta)^{1 - \delta}} \right]^{\mu_{\min}/R}
	\quad\text{for $\delta \in [0, 1]$, and} \\
\Prob{ \lambda_{\max}\left( \sum\nolimits_k \mtx{X}_k \right) \geq (1 + \delta) \mu_{\max} }
	&\leq d \cdot \left[ \frac{\econst^{\delta}}{(1+\delta)^{1 + \delta}} \right]^{\mu_{\max}/R}
	\quad\text{for $\delta \geq 0$.}
\end{align*}
\end{thm}

Chernoff bounds are well suited to studying the spectrum of a random matrix with independent columns.  For additional details and related inequalities, turn to~\S\ref{sec:chernoff}.

\subsubsection{Main Results for Self-Adjoint Matrices}

Another basic example of concentration is provided by a sum of real numbers modulated by independent standard normal variables or, alternatively, by independent Rademacher%
\footnote{A Rademacher random variable is uniformly distributed on $\{\pm 1\}$.}
random variables.  A classical result 
shows that this type of random series exhibits subgaussian tails.  When we replace the real numbers by self-adjoint random matrices, we discover that the maximum and minimum eigenvalue of the matrix sum retain this normal tail behavior.

\begin{thm}[Matrix Gaussian and Rademacher Series] \label{thm:intro-gauss}
Consider a finite sequence $\{ \mtx{A}_k \}$ of fixed, self-adjoint matrices with dimension $d$, and let $\{ \xi_k \}$ be a finite sequence of independent standard normal or independent Rademacher random variables.  Then, for all $t \geq 0$,
$$
\Prob{ \lambda_{\max}\left( \sum\nolimits_k \xi_k \mtx{A}_k \right) \geq t }
	\leq d \cdot \econst^{-t^2/2\sigma^2}
	\quad\text{where}\quad
	\sigma^2 := \norm{ \sum\nolimits_k \mtx{A}_k^2 }.
$$
\end{thm}

Theorem~\ref{thm:intro-gauss} was first established explicitly by Oliveira using a different method~\cite{Oli10:Sums-Random}.  We have included the result here because it is very important and because it follows from a mechanical application of our techniques.  Turn to \S\ref{sec:rademacher} for an exhaustive discussion of matrix Gaussian series.  This presentation also describes several new phenomena that arise when we translate scalar inequalities to the matrix setting.


The Hoeffding inequality is a more general result that describes a sum of independent, zero-mean random variables that are subject to upper and lower bounds; it demonstrates that this random series exhibits normal concentration.  We can extend this result to the matrix setting by considering random matrices that satisfy semidefinite upper bounds.  In the matrix case, the maximum and minimum eigenvalues of the sum also have subgaussian behavior.

\begin{thm}[Matrix Hoeffding] \label{thm:intro-hoeffding}
Consider a finite sequence $\{ \mtx{X}_k \}$ of independent, random, self-adjoint matrices with dimension $d$, and let $\{ \mtx{A}_k \}$ be a sequence of fixed self-adjoint matrices.  Assume that each random matrix satisfies
$$
\Expect \mtx{X}_k = \mtx{0}
\quad\text{and}\quad
\mtx{X}_k^2 \psdle \mtx{A}_k^2
\quad\text{almost surely}.
$$
Then, for all $t \geq 0$,
$$
\Prob{ \lambda_{\max}\left( \sum\nolimits_k \mtx{X}_k \right) \geq t }
	\leq d \cdot \econst^{-t^2 / 8\sigma^2 }
	\quad\text{where}\quad
	\sigma^2 := \norm{ \sum\nolimits_k \mtx{A}_k^2 }.
$$
\end{thm}

The constant $1/8$ in Theorem~\ref{thm:intro-hoeffding} can be improved when there is additional information available.  See~\S\ref{sec:azuma} for a discussion and some related results for martingales.

In fact, a sum of independent, bounded random variables may vary substantially
less than the Hoeffding bound suggests.  A famous inequality of Bernstein
demonstrates that this type of random series exhibits normal concentration near its mean on a scale determined by the
variance of the sum.  On the other hand, the tail of the sum decays subexponentially on a scale
controlled by a uniform upper bound on the summands. Sums of independent random matrices
exhibit the same type of behavior, where the normal concentration depends on a matrix generalization
of the variance and the tails are controlled by a uniform bound on the maximum eigenvalue of each
summand.

\begin{thm}[Matrix Bernstein] \label{thm:intro-bernstein}
Consider a finite sequence $\{ \mtx{X}_k \}$ of independent, random, self-adjoint matrices with dimension $d$.  Assume that each random matrix satisfies
$$
\Expect \mtx{X}_k = \mtx{0}
\quad\text{and}\quad
\lambda_{\max}( \mtx{X}_k ) \leq R
\quad\text{almost surely}.
$$
Then, for all $t \geq 0$,
$$
\Prob{ \lambda_{\max}\left( \sum\nolimits_k \mtx{X}_k \right) \geq t }
	\leq d \cdot \exp\left( \frac{-t^2/2}{\sigma^2 + Rt/3} \right)
	\quad\text{where}\quad
	\sigma^2 := \norm{ \sum\nolimits_k \Expect \big(\mtx{X}_k^2 \big) }.
$$
\end{thm}

Independently, Oliveira has established a somewhat weaker version of Theorem~\ref{thm:intro-bernstein} using alternative techniques~\cite{Oli10:Concentration-Adjacency}.  The reader is probably aware that the probability literature contains a huge number of results that extend Bernstein's inequality to include other \lang{a priori} information on the summands, such as bounds on the rate of moment growth.  Section~\ref{sec:bennett} contains additional matrix probability inequalities of this species.


\subsubsection{Main Results for Rectangular Matrices}

As an immediate corollary of our results for self-adjoint random matrices, we can also establish a collection of inequalities for the maximum singular value of a sum of random rectangular matrices.  In each case, we extend the result to rectangular matrices by using a device from operator theory called the \term{self-adjoint dilation} (\S\ref{sec:dilation}).  Remark~\ref{rem:max-sing} and~\S\ref{sec:gauss-proof} offer some discussion of this technique.  This section presents two of the most important inequalities for sums of random rectangular matrices.

As in the self-adjoint case, the norm of a Gaussian or Rademacher series with rectangular matrix coefficients has subgaussian tails.  This result follows directly from Theorem~\ref{thm:intro-gauss}; see \S\ref{sec:gauss-proof} for a complete proof.  Observe that the variance parameter changes to reflect the fact that the row and column spaces of a general matrix are independent from each other; the variance can be viewed as a noncommutative ``sum of squares.'' 

\begin{thm}[Matrix Gaussian and Rademacher Series: Rectangular Case]
Consider a finite sequence $\{ \mtx{B}_k \}$ of fixed matrices with dimension $d_1 \times d_2$, and let $\{ \xi_k \}$ be a finite sequence of independent standard normal or independent Rademacher random variables.  Define the variance parameter
$$
\sigma^2 := \max\left\{
	\norm{ \sum\nolimits_k \mtx{B}_k \mtx{B}_k^\adj }, \
	\norm{ \sum\nolimits_k \mtx{B}_k^\adj \mtx{B}_k }
\right\}.
$$
Then, for all $t \geq 0$,
$$
\Prob{ \norm{ \sum\nolimits_k \xi_k \mtx{B}_k } \geq t }
	\leq (d_1 + d_2) \cdot \econst^{-t^2/2\sigma^2}.
$$
\end{thm}

We can also develop a rectangular version of the matrix Bernstein inequality.  Notice the parallel between the variance parameter here and the variance parameter for a rectangular Gaussian series.  This result is an immediate corollary of Theorem~\ref{thm:intro-bernstein}; a proof sketch appears in Remark~\ref{rem:rect-bernstein}.

\begin{thm}[Matrix Bernstein: Rectangular Case] \label{thm:intro-bernstein-rect}
Consider a finite sequence $\{ \mtx{Z}_k \}$ of independent, random matrices with dimensions $d_1 \times d_2$.  Assume that each random matrix satisfies
$$
\Expect \mtx{Z}_k = \mtx{0}
\quad\text{and}\quad
\norm{ \mtx{Z}_k } \leq R
\quad\text{almost surely}.
$$
Define
$$
\sigma^2 := \max\left\{ 
	\norm{ \sum\nolimits_k \Expect( \mtx{Z}_k \mtx{Z}_k^\adj ) }, \
	\norm{ \sum\nolimits_k \Expect(\mtx{Z}_k^\adj \mtx{Z}_k) }
\right\}.
$$
Then, for all $t \geq 0$,
$$
\Prob{ \norm{ \sum\nolimits_k \mtx{Z}_k } \geq t }
	\leq (d_1 + d_2) \cdot \exp\left( \frac{-t^2/2}{\sigma^2 + Rt/3} \right).
$$
\end{thm}

We trust that the reader can develop other probability inequalities for rectangular matrices as needed.    For brevity, we have omitted further examples.

\subsubsection{Inequalities for Matrix Martingales}

The techniques in this paper also lead directly to some simple results for matrix martingales.  This material appears in~\S\ref{sec:azuma}.  

\begin{description} \setlength{\itemsep}{0.5pc}
\item[Azuma Inequality]
The Azuma inequality is the martingale extension of the Hoeffding inequality.

\item[McDiarmid Inequality]
The McDiarmid bounded difference inequality concerns matrix-valued functions of a family of independent random variables.  It demonstrates that the extreme eigenvalues of the matrix-valued function exhibit normal concentration.
\end{description}

\noindent
For more refined martingale inequalities, see the papers~\cite{Oli10:Concentration-Adjacency,Tro11:Freedmans-Inequality} and the technical report~\cite{Tro11:User-Friendly-Martingale-TR}.

\subsection{Summary of Related Work} \label{sec:related}



We continue with an overview of some related work on finite-dimensional random matrices.  The first group of papers relies on the matrix extension of the Laplace transform method; the second group uses noncommutative moment inequalities.  

\subsubsection{The Matrix Laplace Transform Method}

The most important precedent for our work is the influential paper of Ahlswede and Winter~\cite{AW02:Strong-Converse}.  They are responsible for developing the matrix version of the Laplace transform method, which shows that the tail probability~\eqref{eqn:intro-prob} is controlled by a matrix generalization of the mgf.  They describe an iterative argument, based on the Golden--Thompson inequality,~\eqref{eqn:golden-thompson} below, that allows them to provide a weak bound for the mgf of a sum of independent random matrices in terms of mgf bounds for the individual summands.  In particular, they apply this technique to obtain an extension of the Chernoff inequality~\cite[Thm.~19]{AW02:Strong-Converse}.


The Ahlswede--Winter method for bounding the matrix mgf is quite general.  Several other authors have exploited their technique to obtain matrix extensions of classical probability inequalities.  Christofides and Markstr{\"o}m establish a matrix version of the Azuma and Hoeffding inequalities~\cite{CM08:Expansion-Properties}.  Gross~\cite[Thm.~6]{Gro11:Recovering-Low-Rank} and Recht~\cite[Thm.~3.2]{Rec09:Simpler-Approach} develop two different matrix extensions of Bernstein's inequality.  We also refer the reader to Vershynin's note~\cite{Ver09:Note-Sums}, which offers a self-contained introduction to the Ahlswede--Winter circle of ideas.

Results established within the Ahlswede--Winter framework are often sharp for sums of i.i.d.~random matrices, but the inequalities are far less accurate when applied to other types of sums.  Roughly speaking, the tail bounds have the correct shape, but the method often leads to poor estimates for the quantity that controls the scale of large deviations.  For a specific example, compare the variance parameter in~\eqref{eqn:aw-intro} with the (correct) variance parameter appearing in~\eqref{eqn:gauss-intro}.  All the results we have mentioned so far have this shortcoming.  See~\S\ref{sec:aw} for technical details.

Very recently, Oliveira has developed two notable variations~\cite{Oli10:Sums-Random,Oli10:Concentration-Adjacency} on the Ahlswede--Winter method for bounding the matrix mgf.  These techniques can sometimes identify the correct matrix generalization of the scale parameter.  In particular, the approach in~\cite{Oli10:Sums-Random} can be used to prove Theorem~\ref{thm:intro-gauss}.  Oliveira has also developed a version of the matrix Bernstein inequality~\cite[Thm.~1.2]{Oli10:Concentration-Adjacency} that is similar to Theorem~\ref{thm:intro-bernstein}; his proof involves a matrix extension of the martingale techniques from~\cite{Fre75:Tail-Probabilities}.

The current article was inspired by the work of Ahlswede--Winter~\cite{AW02:Strong-Converse} and Oliveira~\cite{Oli10:Sums-Random}.  Our results were obtained independently from Oliveira's paper~\cite{Oli10:Concentration-Adjacency}.

\subsubsection{Noncommutative Moment Inequalities}

There is another contemporary line of research that uses noncommutative (nc) moment inequalities to study random matrices.  In a significant article~\cite{Rud99:Random-Vectors}, Rudelson obtains an optimal estimate for the sample complexity of approximating the covariance matrix of a general isotropic distribution.  The argument in his paper, which is due to Pisier, depends on a version of the nc Khintchine inequality~\cite{L-P86:Inegalites-Khintchine,LPP91:Noncommutative-Khintchine,Pis03:Introduction-Operator}.  

Rudelson's technique has been applied widely over the last ten years, and it has emerged as a valuable tool for studying discrete random matrices.  For example, the method can be used to provide bounds on the norm of a random submatrix~\cite[Thm.~1.8]{RV07:Sampling-Large} drawn from a fixed matrix.  It seems likely, however, that matrix probability inequalities will replace the nc Khintchine inequality for many applications because they are easier to use and often produce better results.  

By now, there is a substantial literature on other nc moment inequalities.  The article~\cite{JX05:Best-Constants} contains a reasonably accessible and comprehensive discussion.  Some of these results have been applied to the study of random matrices; see~\cite{JX08:Noncommutative-Burkholder-II} for an example.  As we discuss in~\S\ref{sec:nc-moments}, nc moment bounds can also be combined with the matrix Laplace transform method because they sometimes provide an alternative way to control the matrix mgf.

\subsection{Roadmap}

The rest of the paper is organized as follows.
Section~\ref{sec:prelim} introduces the background results required for our proofs.
Section~\ref{sec:laplace-transform} proves the main technical results that lead to probability inequalities for sums of independent random matrices.  Section~\ref{sec:rademacher} uses Gaussian series as a case study to illustrate the main features of matrix probability inequalities and to argue that the bounds in this paper are structurally optimal.  We develop the matrix Chernoff and Bernstein inequalities in \S\S\ref{sec:chernoff}--\ref{sec:bennett}.  Finally, we establish some simple martingale results in~\S\ref{sec:azuma}.


\section{Algebra, Analysis, and Probability with Matrices} \label{sec:prelim}



This section provides a short introduction to the background we require for our proofs.  The proofs contain detailed cross-references to this material, so the reader may wish to proceed directly to the main thread of argument in~\S\ref{sec:laplace-transform}.





Most of these results can be located in Bhatia's books on matrix analysis~\cite{Bha97:Matrix-Analysis,Bha07:Positive-Definite}.  The works of Horn and Johnson~\cite{HJ85:Matrix-Analysis,HJ94:Topics-Matrix} also serve as good general references.  Higham's book~\cite{Hig08:Functions-Matrices} is an excellent source for information about matrix functions.  


\subsection{Conventions on Matrices}

A \term{matrix} is a finite, two-dimensional array of complex numbers.  \emph{In this paper, all matrices are square unless otherwise noted}.  We add the qualification \term{rectangular} when we need to refer to a general array, which may be square or nonsquare.  Many parts of the discussion do not depend on the size of a matrix, so we specify dimensions only when it matters.  In particular, we usually do not state the size of a matrix when it is determined by the context.

Several abbreviations are ubiquitous.  Instead of self-adjoint, we often write \term{s.a.}  Positive semidefinite becomes \term{psd}, and we shorten positive definite to \term{pd}.

We write $\mtx{0}$ for the zero matrix and $\Id$ for the identity matrix.
The matrix ${\bf E }_{ij}$ has a unit entry in the $(i,j)$ position and zeros elsewhere.  The symbol $\mtx{Q}$ is reserved for a unitary matrix.  We adopt Parlett's convention~\cite{Par87:Symmetric-Eigenvalue} that bold capital letters symmetric about the vertical axis ($\mtx{A}, \dots, \mtx{Y}$ and $\mtx{\Delta}, \dots, \mtx{\Omega}$) refer to s.a.~matrices.  

The symbols $\lambda_{\min}$ and $\lambda_{\max}$ refer to the algebraic minimum and maximum eigenvalues of a s.a.~matrix.  We use curly inequalities to denote the semidefinite ordering: $\mtx{A} \psdge \mtx{0}$ means that $\mtx{A}$ is psd.  The symbol~$\norm{\cdot}$ always refers to the $\ell_2$ vector norm or the associated operator norm, which is called the \term{spectral norm} because it returns the maximum singular value of its argument.



\subsection{Conventions on Probability}

We prefer to avoid unnecessary abstraction and technical detail, so we frame the standing assumption that all random variables are sufficiently regular that we are justified in computing expectations, interchanging limits, and so forth.  Furthermore, we often state that a random variable satisfies some relation and omit the qualification ``almost surely.''  We reserve the symbols $\mtx{X}, \mtx{Y}$ for random s.a.~matrices.

\subsection{Matrix Functions}

Consider a function $f : \mathbb{R} \to \mathbb{R}$.
We define a map on diagonal matrices by applying the function to each diagonal entry.
We then extend $f$ to a function on s.a.~matrices using the eigenvalue decomposition:
\begin{equation} \label{eqn:matrix-fn}
f(\mtx{A}) 
	:= \mtx{Q} \cdot f(\mtx{\Lambda}) \cdot \mtx{Q}^\adj
	\quad\text{where $\mtx{A} = \mtx{Q\Lambda Q}^\adj$.}
\end{equation}
The \term{spectral mapping theorem} states that each eigenvalue of $f(\mtx{A})$ is equal to $f(\lambda)$ for some eigenvalue $\lambda$ of $\mtx{A}$.  This point is obvious from our definition.


Standard inequalities for real functions typically \emph{do not} have parallel versions that hold for the semidefinite ordering.  Nevertheless, there is one type of relation for real functions that always extends to the semidefinite setting:
\begin{equation} \label{eqn:scalar-matrix}
f(a) \leq g(a)
\quad\text{for $a \in I$}
\quad\Longrightarrow\quad
f(\mtx{A}) \psdle g(\mtx{A})
\quad\text{when the eigenvalues of $\mtx{A}$ lie in $I$.}
\end{equation}
We sometimes refer to~\eqref{eqn:scalar-matrix} as the \term{transfer rule}.

\subsection{The Matrix Exponential} \label{sec:matrix-exp}

The exponential of an s.a.~matrix $\mtx{A}$ can be defined by applying~\eqref{eqn:matrix-fn} with the function $f(x) = \econst^x$.  Alternatively, we may use the power series expansion
$$
\exp(\mtx{A}) := 
	\Id + \sum\nolimits_{p=1}^\infty \frac{\mtx{A}^p}{p!}.
$$
The exponential of an s.a.~matrix is always pd~because of the spectral mapping theorem.
On account of the transfer rule~\eqref{eqn:scalar-matrix}, the matrix exponential satisfies some simple semidefinite relations that we collect here.  For each s.a.~matrix $\mtx{A}$, it holds that
\begin{align}
\Id + \mtx{A} &\psdle \econst^{\mtx{A}},
\quad\text{and} \label{eqn:Id+A} \\
\cosh(\mtx{A}) &\psdle \econst^{\mtx{A}^2 / 2}. \label{eqn:cosh-exp}
\end{align}



We often work with the trace of the matrix exponential,
$\trace \exp : \mtx{A} \mapsto \trace \econst^{\mtx{A}}$.
The trace exponential function is convex.  It is also monotone with respect to the semidefinite order:
\begin{equation} \label{eqn:exp-trace-monotone}
\mtx{A} \psdle \mtx{H}
\quad\Longrightarrow\quad
\trace \econst^{\mtx{A}} \leq \trace \econst^{\mtx{H}}.
\end{equation}
See~\cite[Sec.~2]{Pet94:Survey-Certain} for short proofs of these facts.

The matrix exponential \emph{does not} convert sums into products, 
but the trace exponential has a related property that serves as a limited substitute.  The Golden--Thompson inequality~\cite[Sec.~IX.3]{Bha97:Matrix-Analysis} states that
\begin{equation} \label{eqn:golden-thompson}
\trace \econst^{\mtx{A} + \mtx{H}} \leq \trace\left( \econst^{\mtx{A}} \econst^{\mtx{H}} \right)
\quad\text{for all s.a.~$\mtx{A}, \mtx{H}$}.
\end{equation}
The obvious generalization of the bound~\eqref{eqn:golden-thompson} to three matrices is false~\cite[Prob.~IX.8.4]{Bha97:Matrix-Analysis}.

\subsection{The Matrix Logarithm}

We define the matrix logarithm as the functional inverse of the matrix exponential:
\begin{equation} \label{eqn:log-defn}
\log( \econst^{\mtx{A}} ) := \mtx{A}
\quad\text{for each s.a.~matrix $\mtx{A}$}.
\end{equation}
This formula determines the logarithm on the pd cone, which is adequate for our purposes.  

The matrix logarithm interacts beautifully with the semidefinite order~\cite[Exer.~4.2.5]{Bha07:Positive-Definite}.  Indeed, the logarithm is operator monotone:
\begin{equation} \label{eqn:log-monotone}
\mtx{0} \psdlt \mtx{A} \psdle \mtx{H}
\quad\Longrightarrow\quad
\log(\mtx{A}) \psdle \log(\mtx{H}).
\end{equation}
The logarithm is also operator concave:
\begin{equation} \label{eqn:log-concave}
\tau \log(\mtx{A}) + (1-\tau) \log(\mtx{H})
	\psdle \log( \tau \mtx{A} + (1-\tau) \mtx{H} )
\quad\text{for all pd $\mtx{A}, \mtx{H}$ and $\tau \in [0,1]$}.
\end{equation}
{\bf Caveat lector:} Operator monotone functions and operator convex functions are depressingly rare.  In particular, the matrix exponential does not belong to either class~\cite[Ch.~V]{Bha97:Matrix-Analysis}.

\subsection{Dilations} \label{sec:dilation}

An extraordinarily fruitful idea from operator theory is to embed matrices within larger block matrices, called \term{dilations}~\cite{Pau02:Completely-Bounded}.  
The \term{s.a.~dilation} of a rectangular matrix $\mtx{B}$ is
\begin{equation} \label{eqn:sa-dilation}
\coll{S}(\mtx{B}) := \begin{bmatrix} \mtx{0} & \mtx{B} \\ \mtx{B}^\adj & \mtx{0} \end{bmatrix}.
\end{equation}
Evidently, $\coll{S}(\mtx{B})$ is always s.a.  
A short calculation yields the important identity
\begin{equation} \label{eqn:sa-modulus}
\coll{S}( \mtx{B} )^2
	= \begin{bmatrix} \mtx{BB}^\adj & \mtx{0} \\
	\mtx{0} & \mtx{B}^\adj \mtx{B} \end{bmatrix}.
\end{equation}
It can also be verified that the s.a.~dilation preserves spectral information:
\begin{equation} \label{eqn:sa-norm}
\lambda_{\max}(\coll{S}(\mtx{B})) = \norm{ \coll{S}(\mtx{B}) } = \norm{ \mtx{B} }.
\end{equation}
We use dilations to extend results for s.a.~matrices to rectangular matrices.  See Remark~\ref{rem:max-sing} and~\S\ref{sec:gauss-proof} for more information about this technique.

\subsection{Expectation and the Semidefinite Order}

Since the expectation of a random matrix can be viewed as a convex combination and the psd~cone is convex, expectation preserves the semidefinite order:
\begin{equation} \label{eqn:expect-psd}
\mtx{X} \psdle \mtx{Y}
\quad\text{almost surely}
\quad\Longrightarrow\quad
\Expect \mtx{X} \psdle \Expect \mtx{Y}.
\end{equation}



Every operator convex function admits an operator Jensen's inequality~\cite{HP03:Jensens-Operator}.  In particular, the matrix square is operator convex, which implies that
\begin{equation} \label{eqn:jensen-matrix-square}
(\Expect \mtx{X})^2 \psdle \Expect \bigl(\mtx{X}^2 \bigr).
\end{equation}
The relation~\eqref{eqn:jensen-matrix-square} is also a specific instance of Kadison's inequality~\cite[Thm.~2.3.2]{Bha07:Positive-Definite}.

\section{Tail Bounds via the Laplace Transform Method} \label{sec:laplace-transform}

This section develops some general probability inequalities for the maximum eigenvalue of a sum of independent random matrices.  The main argument can be viewed as a matrix extension of the Laplace transform method for sums of independent real random variables.  In the matrix setting, however, it requires great care to execute this technique successfully.

\subsection{Matrix Moments and Cumulants} \label{sec:mom-cum}

Consider a random s.a.~matrix $\mtx{X}$ that has moments of all orders.
By analogy with the classical scalar definitions, we may construct matrix extensions of the moment generating function (mgf) and the cumulant generating function (cgf):
\begin{equation} \label{eqn:matrix-mgf-cgf}
\mtx{M}_{\mtx{X}}(\theta) := \Expect \econst^{\theta \mtx{X}}
\quad\text{and}\quad
\mtx{\Xi}_{\mtx{X}}(\theta) := \log \Expect \econst^{\theta \mtx{X}}
\quad\text{for $\theta \in \mathbb{R}$.}
\end{equation}
We admit the possibility that these expectations do not exist for all values of $\theta$.  The matrix cgf can be viewed as an \term{exponential mean}, a weighted average that emphasizes large deviations (with the same sign as $\theta$).
The matrix mgf and cgf have formal power series expansions:
$$
\mtx{M}_{\mtx{X}}(\theta) = \Id + \sum\nolimits_{p=1}^\infty
	\frac{\theta^p }{ p! } \cdot\Expect( \mtx{X}^p )
\quad\text{and}\quad
\mtx{\Xi}_{\mtx{X}}(\theta) = \sum\nolimits_{p=1}^\infty \frac{\theta^p}{p!} \cdot \mtx{\Psi}_p.
$$
The coefficients $\Expect (\mtx{X}^p)$ are called \term{matrix moments}, and we refer to $\mtx{\Psi}_p$ as a \term{matrix cumulant}.
The matrix cumulant $\mtx{\Psi}_p$ has a formal expression as a (noncommutative) polynomial in the matrix moments up to order $p$.  In particular, the first cumulant is the mean and the second cumulant is the variance:
$$
\mtx{\Psi}_1 = \Expect \mtx{X}
\quad\text{and}\quad
\mtx{\Psi}_2 = \Expect (\mtx{X}^2) - (\Expect \mtx{X})^2.
$$
Higher-order cumulants are harder to write down and interpret.


\subsection{The Laplace Transform Method for Matrices} \label{sec:matrix-lt}


We begin our main development with a striking idea drawn from the influential paper~\cite{AW02:Strong-Converse} of Ahlswede and Winter.  Their work contains a matrix analog of the classical Laplace transform bound. We need the following variant, which is due to Oliveira~\cite{Oli10:Sums-Random}.

\begin{prop}[The Laplace Transform Method] \label{prop:laplace-transform}
Let $\mtx{Y}$ be a random self-adjoint~matrix.  For all $t \in \mathbb{R}$,
$$
\Prob{ \lambda_{\max}(\mtx{Y}) \geq t }
	\leq \inf_{\theta > 0} \left\{ \econst^{-\theta t}
	\cdot \Expect \trace \econst^{\theta \mtx{Y}} \right\}.
$$
\end{prop}

In words, we can control tail probabilities for the maximum eigenvalue of a random matrix by producing a bound for the trace of the matrix mgf defined in~\eqref{eqn:matrix-mgf-cgf}.  

\begin{proof}
Fix a positive number $\theta$.  We have the chain of relations
$$
\Prob{ \lambda_{\max}(\mtx{Y}) \geq t }
	= \Prob{ \lambda_{\max}(\theta \mtx{Y}) \geq \theta t }
	= \Prob{ \econst^{\lambda_{\max}(\theta \mtx{Y})} \geq \econst^{\theta t} }
	\leq \econst^{- \theta t} \cdot \Expect \econst^{\lambda_{\max}(\theta\mtx{Y})}.
$$
The first identity uses the homogeneity of the maximum eigenvalue map, and the second relies on the monotonicity of the scalar exponential function; the third relation is Markov's inequality.  To bound the exponential, note that
$$
\econst^{\lambda_{\max}(\theta \mtx{Y})}
	= \lambda_{\max}\bigl(\econst^{\theta\mtx{Y}}\bigr)
	\leq \trace \econst^{\theta\mtx{Y}}.
$$
The identity is the spectral mapping theorem; the inequality holds because the exponential of an s.a.~matrix is pd~and the maximum eigenvalue of a pd~matrix is dominated by the trace.  Combine the latter two relations to reach
$$
\Prob{ \lambda_{\max}(\mtx{Y}) \geq t }
	\leq \econst^{- \theta t} \cdot \Expect \trace \econst^{\theta\mtx{Y}}.
$$
This inequality holds for any positive $\theta$, so we may take an infimum to complete the proof.
\end{proof}


\subsection{The Failure of the Matrix mgf} \label{sec:mom-fail}

%


In the scalar setting, the Laplace transform method is very effective for studying sums of independent random variables because the mgf decomposes.
Consider an independent sequence $\{X_k\}$ of real random variables.  Operating formally, we see that the (scalar) mgf of the sum satisfies a multiplication rule:
\begin{equation} \label{eqn:mgf-mult}
M_{(\sum\nolimits_k X_k)}(\theta)
	= \Expect \exp\left(\sum\nolimits_k \theta X_k \right)
	= \Expect \prod\nolimits_k \econst^{\theta X_k}
	= \prod\nolimits_k \Expect \econst^{\theta X_k}
	= \prod\nolimits_k M_{X_k}(\theta).
\end{equation}
{\em This calculation relies on the fact that the scalar exponential function converts sums to products, a property the matrix exponential does not share.}  As a consequence, there is no immediate analog of~\eqref{eqn:mgf-mult} in the matrix setting.  


Ahlswede and Winter attempt to imitate the multiplication rule~\eqref{eqn:mgf-mult} using the following observation.  When $\mtx{X}_1$ and $\mtx{X}_2$ are independent random matrices,
\begin{equation} \label{eqn:aw-mgf}
\trace \mtx{M}_{\mtx{X}_1 + \mtx{X}_2}(\theta)
	\leq \Expect \trace \big[ \econst^{\theta \mtx{X}_1} \econst^{\theta \mtx{X}_2} \big]
	= \trace\big[ (\Expect \econst^{\theta \mtx{X}_1})(\Expect \econst^{\theta \mtx{X}_2})\big]
	= \trace \big[ \mtx{M}_{\mtx{X}_1}(\theta) \cdot \mtx{M}_{\mtx{X}_2}(\theta) \big].
\end{equation}
The first relation is the Golden--Thompson trace inequality~\eqref{eqn:golden-thompson}.  Unfortunately, we cannot extend the bound~\eqref{eqn:aw-mgf} to include additional matrices.
This cold fact suggests that the Golden--Thompson inequality may not be the natural way to proceed.  In~\S\ref{sec:aw}, we map out the route Ahlswede and Winter pursue, but we continue along a different path.

\subsection{A Concave Trace Function} \label{sec:lieb}



For inspiration, we turn to the literature on matrix analysis.  Some of the most beautiful and profound results in this domain concern the convexity of trace functions.  We have observed that this theory has incredible implications for the study of random matrices.  This paper demonstrates that a large class of matrix probability inequalities follows from a deep theorem~\cite[Thm.~6]{Lie73:Convex-Trace} of Lieb that appears in his seminal work on convex trace functions.

\begin{thm}[Lieb] \label{thm:lieb}
Fix a self-adjoint~matrix $\mtx{H}$.  The function
$$
\mtx{A} \longmapsto \trace \exp( \mtx{H} + \log(\mtx{A}))
$$
is concave on the positive-definite~cone.
\end{thm}

Epstein provides an alternative proof of Theorem~\ref{thm:lieb} in~\cite[Sec.~II]{Eps73:Remarks-Two}, and Ruskai offers a simplified account of Epstein's argument in~\cite{Rus02:Inequalities-Quantum,Rus05:Erratum-Inequalities}.  The note~\cite{Tro11:Joint-Convexity} derives Lieb's theorem from the joint convexity of quantum relative entropy~\cite[Lem.~2]{Lin74:Expectations-Entropy}. The latter approach is advantageous because the joint convexity result admits several elegant, conceptual proofs, such as~\cite[Cor.~2.2]{Eff09:Matrix-Convexity}.




We require a simple but powerful corollary of Lieb's theorem.  This result describes how expectation interacts with the trace exponential.


\begin{cor} 
\label{cor:cum-ineq}
Let $\mtx{H}$ be a fixed self-adjoint matrix, and let $\mtx{X}$ be a random self-adjoint matrix.  Then
$$
\Expect \trace \exp( \mtx{H} + \mtx{X} )
	\leq \trace \exp( \mtx{H} + \log( \Expect \econst^{\mtx{X}} ) ).
$$
\end{cor}

\begin{proof}
Define the random matrix $\mtx{Y} = \econst^{\mtx{X}}$, and calculate that
$$
\Expect \trace \exp(\mtx{H} + \mtx{X})
	= \Expect \trace \exp( \mtx{H} + \log( \mtx{Y} ) )
	\leq \trace \exp( \mtx{H} + \log( \Expect \mtx{Y} ) )
	= \trace \exp( \mtx{H} + \log( \Expect \econst^{\mtx{X}} ) ).
$$
The first identity follows from the definition~\eqref{eqn:log-defn} of the matrix logarithm because $\mtx{Y}$ is always pd.  Lieb's result, Theorem~\ref{thm:lieb}, ensures that the trace function is concave in $\mtx{Y}$, so we may invoke Jensen's inequality to draw the expectation inside the logarithm.
\end{proof}

\subsection{Subadditivity of the Matrix cgf} \label{sec:cum-subadd}

Let us return to the problem of bounding the matrix mgf of an independent sum.
Although the multiplication rule~\eqref{eqn:mgf-mult} is a dead end in the matrix case, the scalar cgf has a related property that submits to generalization.  For an independent family $\{X_k\}$ of real random variables, the scalar cgf is additive:
\begin{equation} \label{eqn:cgf-add}
\Xi_{(\sum\nolimits_k X_k)}(\theta)
	= \log \Expect \exp\left( \sum\nolimits_k \theta X_k \right)
	= \sum\nolimits_k \log \Expect \econst^{\theta X_k}
	= \sum\nolimits_k \Xi_{X_k}(\theta),
\end{equation}
where the second identity follows from~\eqref{eqn:mgf-mult} when we take logarithms.

Our key insight is that Corollary~\ref{cor:cum-ineq} offers a completely satisfactory way to extend the addition rule~\eqref{eqn:cgf-add} for scalar cgfs to the matrix setting.  We have the following result.

%
%
%



\begin{lemma}[Subadditivity of Matrix cgfs] \label{lem:cgf-indep}
Consider a finite sequence $\{ \mtx{X}_k \}$ of independent, random, self-adjoint matrices.  Then
$$
\Expect \trace \exp\left( \sum\nolimits_k \theta \mtx{X}_k \right)
	\leq \trace \exp\left( \sum\nolimits_k \log \Expect \econst^{\theta \mtx{X}_k} \right)
	\quad\text{for $\theta \in \mathbb{R}$.}
$$
\end{lemma}


\begin{proof}
It does no harm to assume $\theta = 1$.
Let $\Expect_k$ denote the expectation, conditioned on $\mtx{X}_1, \dots, \mtx{X}_{k}$.  Abbreviate
\begin{equation*}\label{eqn:xi-k}
\mtx{\Xi}_k := \log( \Expect_{k-1} \econst^{\mtx{X}_k} )
	= \log( \Expect \econst^{\mtx{X}_k} ),
\end{equation*}
where the equality holds because the family $\{\mtx{X}_k\}$ is independent.   We see that
\begin{align*}
\Expect \trace \exp\left( \sum\nolimits_{k=1}^n \mtx{X}_k \right)
	&= \Expect_0 \cdots \Expect_{n-1}
		\trace \exp\left( \sum\nolimits_{k=1}^{n-1} \mtx{X}_k
		+ \mtx{X}_n \right)  \\
	&\leq \Expect_{0} \cdots \Expect_{n-2}
		\trace \exp\left( \sum\nolimits_{k=1}^{n-1} \mtx{X}_k
		+ \log(\Expect_{n-1} \econst^{ \mtx{X}_n} ) \right) \\
	&= \Expect_{0} \cdots \Expect_{n-2}
		\trace \exp\left( \sum\nolimits_{k=1}^{n-2} \mtx{X}_k
		+ \mtx{X}_{n-1} + \mtx{\Xi}_n \right) \\
	&\leq \Expect_{0} \cdots \Expect_{n-3}
		\trace \exp\left( \sum\nolimits_{k=1}^{n-2} \mtx{X}_k
		+ \mtx{\Xi}_{n-1} + \mtx{\Xi}_{n} \right) \\
\dots\quad
	&\leq \trace \exp\left( \sum\nolimits_{k=1}^{n} \mtx{\Xi}_k \right).
\end{align*}
The first line relies on the tower property of conditional expectation.  At each step $m = 1, 2, \dots, n$, we invoke Corollary~\ref{cor:cum-ineq} with the fixed matrix $\mtx{H}$ equal to
$$
\mtx{H}_m = \sum\nolimits_{k=1}^{m-1} \mtx{X}_k +
	\sum\nolimits_{k=m+1}^n \mtx{\Xi}_k.
$$
This act is legal because $\mtx{H}_m$ does not depend on $\mtx{X}_m$.
\end{proof}

\begin{rem}
To make the parallel with the addition rule~\eqref{eqn:cgf-add} clearer, we can rewrite the conclusion of Lemma~\ref{lem:cgf-indep} in the form
$$
\trace \exp\left( \mtx{\Xi}_{(\sum_k \mtx{X}_k)}(\theta) \right)
	\leq \trace \exp\left( \sum\nolimits_k \mtx{\Xi}_{\mtx{X}_k}(\theta) \right)
$$
by applying the definition~\eqref{eqn:matrix-mgf-cgf} of the matrix cgf.
\end{rem}

\subsection{Tail Bounds for Independent Sums} \label{sec:master-tail}

This section contains abstract tail bounds for the sum of independent random matrices.  Later, we will specialize these results to some specific situations.  We begin with a very general inequality, which is the progenitor of our other results.

\begin{thm}[Master Tail Bound for Independent Sums] \label{thm:master-ineq}
Consider a finite sequence $\{ \mtx{X}_k \}$ of independent, random, self-adjoint matrices.  For all $t \in \mathbb{R}$,
\begin{equation} \label{eqn:master-ineq}
\Prob{ \lambda_{\max}\left( \sum\nolimits_k \mtx{X}_k \right) \geq t }
	\leq \inf_{\theta > 0} \left\{ \econst^{-\theta t}
	\cdot \trace \exp\left( \sum\nolimits_k
	 \log \Expect \econst^{\theta \mtx{X}_k} \right)
	\right\}.
\end{equation}
\end{thm}

\begin{proof}
Substitute the subadditivity rule for matrix cgfs, Lemma~\ref{lem:cgf-indep}, into the Laplace transform bound, Proposition~\ref{prop:laplace-transform}.
\end{proof}

Our first corollary adapts Theorem~\ref{thm:master-ineq} to the case that arises most often in practice.  We call upon this result several times to obtain tail bounds under a variety of assumptions about the structure of the random matrices.

\begin{cor} \label{cor:main-result-indep}
Consider a finite sequence $\{ \mtx{X}_k \}$ of independent, random, self-adjoint~matrices with dimension $d$.   Assume there is a function $g : (0, \infty) \to [0, \infty]$ and a sequence $\{\mtx{A}_k\}$ of fixed self-adjoint matrices that satisfy the relations
\begin{equation} \label{eqn:mgf-hypothesis}
\Expect \econst^{\theta \mtx{X}_k} \psdle \econst^{g(\theta) \cdot \mtx{A}_k}
\quad\text{for $\theta > 0$.}
\end{equation}
Define the scale parameter
$$
\rho := \lambda_{\max}\left( \sum\nolimits_k \mtx{A}_k \right).
$$
Then, for all $t \in \mathbb{R}$,
\begin{equation} \label{eqn:deploy}
\Prob{ \lambda_{\max}\left(\sum\nolimits_k \mtx{X}_k \right) \geq t }
	\leq d \cdot \inf_{\theta > 0} \econst^{- \theta t + g(\theta) \cdot \rho }.
\end{equation}
\end{cor}



\begin{proof}
The hypothesis~\eqref{eqn:mgf-hypothesis} implies that
\begin{equation} \label{eqn:cgf-hypothesis}
\log \Expect \econst^{\theta \mtx{X}_k} \psdle g(\theta) \cdot \mtx{A}_k
\quad\text{for $\theta > 0$}
\end{equation}
because of the property~\eqref{eqn:log-monotone} that the matrix logarithm is operator monotone.
Recall the fact~\eqref{eqn:exp-trace-monotone} that the trace exponential is monotone with respect to the semidefinite order.  As a consequence, we can introduce each relation from the family~\eqref{eqn:cgf-hypothesis} into the master inequality~\eqref{eqn:master-ineq}.  For each $\theta > 0$, it follows that
\begin{align*}
\Prob{ \lambda_{\max}\left(\sum\nolimits_k \mtx{X}_k \right) \geq t }
	&\leq \econst^{-\theta t}
	\cdot \trace \exp\left( g(\theta) \cdot \sum\nolimits_k \mtx{A}_k \right) \\
	&\leq \econst^{-\theta t}
	\cdot d \cdot \lambda_{\max}\left( \exp \left( g(\theta) \cdot \sum\nolimits_k \mtx{A}_k \right) \right) \\
	&= d \cdot \econst^{-\theta t}
	\cdot \exp\left( g(\theta) \cdot \lambda_{\max}\left( \sum\nolimits_k \mtx{A}_k
	\right)\right).	
\end{align*}
The second inequality holds because the trace of a pd matrix, such as the exponential, is bounded by the dimension $d$ times the maximum eigenvalue.  The last line depends on the spectral mapping theorem and the fact that the function $g$ is nonnegative.  Identify the quantity $\rho$, and take the infimum over positive $\theta$ to reach the conclusion~\eqref{eqn:deploy}.
\end{proof}

\begin{rem}
An alternative expression of the result~\eqref{eqn:deploy} is that
\begin{equation*} 
\Prob{ \lambda_{\max}\left(\sum\nolimits_k \mtx{X}_k \right) \geq t }
	\leq d \cdot \exp\left( - \sup_{\theta > 0} \left\{ \theta t - g(\theta) \cdot \rho \right\} \right)
	= d \cdot \exp\left( - \rho \cdot g^*(t/\rho) \right).
\end{equation*}
In words, the exponent in the tail bound can be written in terms of the perspective transformation of the Fenchel--Legendre conjugate of the function $g$. 
This inequality parallels the upper estimate in Cram{\'e}r's classical result for large deviations~\cite[Thm.~2.2.3]{DZ98:Large-Deviations}.  
\end{rem}

It is also worthwhile to state another consequence of Theorem~\ref{thm:master-ineq}.  This bound is sometimes more useful than Corollary~\ref{cor:main-result-indep} because it combines the mgfs of the random matrices together under a single logarithm.

\begin{cor} \label{cor:mgf-master-ineq} 
Consider a sequence $\{ \mtx{X}_k : k = 1, 2, \dots, n \}$ of independent, random, self-adjoint matrices with dimension $d$.  For all $t \in \mathbb{R}$,
\begin{equation} \label{eqn:mgf-master-ineq}
\Prob{ \lambda_{\max}\left( \sum\nolimits_{k=1}^n \mtx{X}_k \right) \geq t }
	\leq d \cdot \inf_{\theta > 0} \exp\left( - \theta t +
	n \cdot \log \lambda_{\max}\left(
	\frac{1}{n} \sum\nolimits_{k=1}^n \Expect \econst^{\theta \mtx{X}_k} \right)
	\right).
\end{equation}
\end{cor}

\begin{proof}
Recall the fact~\eqref{eqn:log-concave} that the matrix logarithm is operator concave.  For each $\theta > 0$, it follows that
$$
\sum\nolimits_{k=1}^n \log \Expect \econst^{\theta \mtx{X}_k}
	= n \cdot \frac{1}{n} \sum\nolimits_{k=1}^n \log \Expect \econst^{\theta \mtx{X}_k}
	\psdle n \cdot \log \left( \frac{1}{n} \sum\nolimits_{k=1}^n \Expect \econst^{\theta \mtx{X}_k} \right).
$$
The property~\eqref{eqn:exp-trace-monotone} that the trace exponential is monotone allows us to introduce the latter relation into the master inequality~\eqref{eqn:master-ineq} to obtain
$$
\Prob{ \lambda_{\max}\left( \sum\nolimits_{k=1}^n \mtx{X}_k \right) \geq t }
	\leq \econst^{-\theta t}
	\cdot \trace \exp \left( n \cdot \log \left( \frac{1}{n} \sum\nolimits_{k=1}^n
	 \Expect \econst^{\theta \mtx{X}_k} \right) \right).
$$
To complete the proof, we bound the trace by $d$ times the maximum eigenvalue, and we invoke the spectral mapping theorem (twice!) to draw the maximum eigenvalue map inside the logarithm.  Take the infimum over positive $\theta$ to reach~\eqref{eqn:mgf-master-ineq}.
\end{proof}

We conclude this section with remarks on some other situations that we can analyze using the master tail bound, Theorem~\ref{thm:master-ineq}, and its corollaries.

\begin{rem}[Minimum Eigenvalue]
We can study the minimum eigenvalue of a sum of random s.a.~matrices because
$
\lambda_{\min}(\mtx{X})
	= - \lambda_{\max}(-\mtx{X}).
$
As a result,
$$
\Prob{ \lambda_{\min}\left(\sum\nolimits_k \mtx{X}_k \right) \leq t }
	= \Prob{ \lambda_{\max}\left(\sum\nolimits_k -\mtx{X}_k \right) \geq -t }.
$$
In \S\ref{sec:chernoff}, we apply this observation to develop lower Chernoff bounds.
\end{rem}

\begin{rem}[Maximum Singular Value] \label{rem:max-sing}
We can also analyze the maximum singular value of a sum of random rectangular matrices by applying these results to the s.a.~dilation~\eqref{eqn:sa-dilation}.  For a finite sequence $\{\mtx{Z}_k\}$ of independent, random, rectangular matrices, we have
$$
\Prob{ \norm{ \sum\nolimits_k \mtx{Z}_k } \geq t }
	= \Prob{ \lambda_{\max}\left(\sum\nolimits_k \coll{S}(\mtx{Z}_k) \right) \geq t }
$$
on account of~\eqref{eqn:sa-norm} and the property that the dilation is real-linear.
This device allows us to extend most of the tail bounds in this paper to rectangular matrices.  See~\S\ref{sec:rademacher} for an application to Gaussian and Rademacher series.
\end{rem}

\begin{rem}[Martingales]
It is possible to combine the proofs of Lemma~\ref{lem:cgf-indep} and Theorem~\ref{thm:master-ineq} to obtain some simple results for matrix martingales.  See the demonstration of the matrix Azuma inequality in~\S\ref{sec:azuma} for an example of this approach.  To reach fully detailed results for martingales, one must use a fundamentally different style of argument~\cite{Oli10:Concentration-Adjacency,Tro11:Freedmans-Inequality}.
\end{rem}

\subsection{The Ahlswede--Winter Method} \label{sec:aw}

Ahlswede and Winter use a different approach to bound the matrix mgf, which exploits the multiplicative bound~\eqref{eqn:aw-mgf} for the trace exponential of a sum of two independent, random, s.a.~matrices.  The reader may find their argument interesting.

Consider a sequence $\{\mtx{X}_k : k = 1, 2, \dots, n \}$ of independent, random, s.a.~matrices with dimension $d$, and let $\mtx{Y} = \sum\nolimits_k \mtx{X}_k$.  The trace inequality~\eqref{eqn:aw-mgf} implies that
$$
\trace \mtx{M}_{\mtx{Y}}(\theta) \leq \trace \left[ \big(\Expect \econst^{\sum_{k=1}^{n-1} \theta \mtx{X}_k} \big)\big(\Expect \econst^{\theta \mtx{X}_n} \big) \right]
	\leq \trace \big(\Expect \econst^{\sum_{k=1}^{n-1} \theta \mtx{X}_k} \big) 
	\cdot \lambda_{\max}\big( \Expect \econst^{\theta \mtx{X}_n} \big).
$$
Iterating this procedure leads to the relation
\begin{equation} \label{eqn:aw-mgf-bd}
\trace \mtx{M}_{\mtx{Y}}(\theta)
	\leq (\trace \Id) \cdot \left[ \prod\nolimits_k \lambda_{\max}\big( \Expect \econst^{\theta \mtx{X}_k} \big)
	\right]
	= d \cdot \exp\left(\sum\nolimits_k \lambda_{\max}\big( \log \Expect  \econst^{\theta \mtx{X}_k} \big) \right).
\end{equation}
The bound~\eqref{eqn:aw-mgf-bd} is the key to the Ahlswede--Winter method for producing probability inequalities.  As a consequence, their approach generally leads to tail bounds that depend on a scale parameter involving ``the sum of eigenvalues.''  See, for example, the bound~\eqref{eqn:aw-intro} or the matrix probability inequalities presented in the papers~\cite{AW02:Strong-Converse,CM08:Expansion-Properties,Gro11:Recovering-Low-Rank,Rec09:Simpler-Approach}.

In contrast, our result on the subadditivity of cumulants, Lemma~\ref{lem:cgf-indep}, implies that
\begin{equation} \label{eqn:my-mgf-bd}
\trace \mtx{M}_{\mtx{Y}}(\theta)
	\leq d \cdot \exp \left(\lambda_{\max}\left( \sum\nolimits_k \log \Expect \econst^{\theta \mtx{X}_k} \right) \right).
\end{equation}
Probability inequalities developed with~\eqref{eqn:my-mgf-bd} contain a scale parameter that involves the ``eigenvalue of a sum.'' See, for example, the bound~\eqref{eqn:gauss-intro}.
The exponent in~\eqref{eqn:aw-mgf-bd} often exceeds the exponent in~\eqref{eqn:my-mgf-bd} by a factor of $d$, the ambient dimension, which is a serious loss.  Section~\ref{sec:aw-gauss} describes concrete situations where this discrepancy occurs.

\section{Case Study: Matrix Gaussian Series} \label{sec:rademacher}



A matrix Gaussian series stands among the simplest instances of a sum of independent random matrices.  Nevertheless, this example already exhibits several new phenomena that arise when we translate scalar tail bounds to the matrix setting.  Consequently, we explore this fundamental case in depth as a way to develop insights about other matrix probability inequalities.

\subsection{Main Results}

We begin with the scalar case.  Consider a finite sequence $\{a_k\}$ of real numbers and a finite sequence $\{\gamma_k\}$ of independent standard Gaussian variables.  We have the probability inequality
\begin{equation} \label{eqn:real-gauss}
\Prob{ \sum\nolimits_{k} \gamma_k \, a_k \geq t }
	\leq \econst^{- t^2 / 2\sigma^2 }
\quad\text{where $\sigma^2 := \sum\nolimits_{k} a_k^2$.}
\end{equation}
This result testifies that a Gaussian series with real coefficients satisfies a normal-type tail bound where the variance is controlled by the sum of the squared coefficients.
The relation~\eqref{eqn:real-gauss} follows easily from the scalar Laplace transform method.  An alternative proof proceeds using the rotational invariance of a standard normal vector along with basic estimates on the error function.

The inequality~\eqref{eqn:real-gauss} generalizes directly to the noncommutative setting, as do many other scalar tail bounds.  The matrix Laplace transform method, Proposition~\ref{prop:laplace-transform}, delivers the following result on the tail behavior of a matrix Gaussian series.

\begin{thm}[Matrix Gaussian and Rademacher Series] \label{thm:rad-gauss-series}
Consider a finite sequence $\{ \mtx{A}_k \}$ of fixed self-adjoint~matrices with dimension $d$, and let $\{\gamma_k\}$ be a finite sequence of independent standard normal variables.  Compute the variance parameter
\begin{equation} \label{eqn:matrix-gauss-sd}
\sigma^2 := \norm{ \sum\nolimits_k \mtx{A}_k^2 }.
\end{equation}
Then, for all $t \geq 0$,
\begin{equation} \label{eqn:matrix-gauss}
\Prob{ \lambda_{\max}\left( \sum\nolimits_k \gamma_k \mtx{A}_k \right) \geq t}
	\leq d \cdot \econst^{- t^2 / 2\sigma^2}.
\end{equation}
In particular,
\begin{equation} \label{eqn:twosided-matrix-gauss}
\Prob{ \norm{ \sum\nolimits_k \gamma_k \mtx{A}_k } \geq t}
	\leq 2d \cdot \econst^{- t^2 / 2\sigma^2}.
\end{equation}
The same bounds hold when we replace $\{\gamma_k\}$ by a finite sequence of independent Rademacher random variables.
\end{thm}

Observe that the bound~\eqref{eqn:matrix-gauss} reduces to the scalar result~\eqref{eqn:real-gauss} when the dimension $d = 1$.
Of course, one may wonder whether the generalization~\eqref{eqn:matrix-gauss-sd} of the scalar variance is sharp and whether the dimensional dependence in~\eqref{eqn:matrix-gauss} is necessary.  A primary objective of this section is to demonstrate that Theorem~\ref{thm:rad-gauss-series} cannot be improved without 
changing its form.



Most of the inequalities in this paper have variants that concern the maximum singular value of a sum of rectangular random matrices.
These extensions follow immediately when we apply the s.a.~results to the s.a.~dilation of the sum of rectangular matrices.  Here is the general version of~Theorem~\ref{thm:rad-gauss-series}, which serves as a model for other rectangular results.


\begin{cor}[Rectangular Matrix Gaussian and Rademacher Series] \label{cor:rect-rad-gauss-series}
Consider a finite sequence $\{ \mtx{B}_k \}$ of fixed matrices with dimension $d_1 \times d_2$, and let $\{\gamma_k\}$ be a finite sequence of independent standard normal variables.  Compute the variance parameter
\begin{equation*} 
\sigma^2 := \max\left\{ \norm{ \sum\nolimits_k \mtx{B}_k \mtx{B}_k^\adj }, \
	\norm{ \sum\nolimits_k \mtx{B}_k^\adj \mtx{B}_k } \right\}.
\end{equation*}
Then, for all $t \geq 0$,
\begin{equation*}
\Prob{ \norm{ \sum\nolimits_k \gamma_k \mtx{B}_k } \geq t}
	\leq (d_1 + d_2) \cdot \econst^{- t^2 / 2\sigma^2}.
\end{equation*}
The same bound holds when we replace $\{\gamma_k\}$ by a finite sequence of independent Rademacher random variables.
\end{cor}

The proofs of Theorem~\ref{thm:rad-gauss-series} and Corollary~\ref{cor:rect-rad-gauss-series} appear below in~\S\ref{sec:gauss-proof}.  Unlike our other results, these two bounds are not new.  One established argument, which we discuss in~\S\ref{sec:nc-moments}, involves noncommutative Khintchine inequalities.  It is also possible to prove these results using Oliveira's ideas~\cite{Oli10:Sums-Random}.



\subsection{Proofs} \label{sec:gauss-proof}

We continue with a short demonstration of the main results for matrix Gaussian and Rademacher series.  The first step is to obtain a semidefinite bound for the mgf of a fixed matrix modulated by a Gaussian variable or a Rademacher variable.  This mgf bound essentially appears in Oliveira's work~\cite[Lem.~2]{Oli10:Sums-Random}.

\begin{lemma}[Rademacher and Gaussian mgfs] \label{lem:rad-gauss-mgf}
Suppose that $\mtx{A}$ is an s.a.~matrix.  Let $\eps$ be a Rademacher random variable, and let $\gamma$ be a standard normal random variable.  Then
$$
\Expect \econst^{\eps \theta \mtx{A}}
\psdle \econst^{\theta^2\mtx{A}^2/2}
\quad\text{and}\quad
\Expect \econst^{\gamma \theta \mtx{A}}
= \econst^{\theta^2 \mtx{A}^2/2}
\quad\text{for $\theta \in \mathbb{R}$.}
$$ 
\end{lemma}

\begin{proof}
Absorbing $\theta$ into $\mtx{A}$, we may assume $\theta = 1$ in each case.  We begin with the Rademacher mgf.  By direct calculation,
$$
\Expect \econst^{\eps \mtx{A}}
	= \cosh(\mtx{A})
	\psdle \econst^{\mtx{A}^2/2},
$$
where the second relation is~\eqref{eqn:cosh-exp}.

For the Gaussian case, recall that the moments of a standard normal variable satisfy
$$
\Expect( \gamma^{2p+1} ) = 0
\quad\text{and}\quad
\Expect( \gamma^{2p} ) = \frac{(2p)!}{p! \, 2^p}
\quad\text{for $p = 0, 1, 2, \dots$}.
$$
Therefore,
$$
\Expect \econst^{\gamma \mtx{A}}
	= \Id + \sum_{p=1}^\infty \frac{\Expect(\gamma^{2p}) \mtx{A}^{2p}}{(2p)!}
	= \Id + \sum_{p=1}^\infty \frac{(\mtx{A}^2/2)^p}{p!}
	= \econst^{ \mtx{A}^2 / 2 }.
$$
The first identity holds because the odd terms in the series vanish.
\end{proof}

The tail bounds for s.a.~matrix Gaussian and Rademacher series follow easily.


\begin{proof}[Proof of Theorem~\ref{thm:rad-gauss-series}]
Let $\{\xi_k\}$ be a finite sequence of independent standard normal variables or independent Rademacher variables.
Invoke Lemma~\ref{lem:rad-gauss-mgf} to obtain 
$$
\Expect \econst^{\xi_k \theta \mtx{A}_k}
\psdle \econst^{g(\theta) \cdot \mtx{A}_k^2}
\quad\text{where $g(\theta) := \theta^2/2$ for $\theta > 0$.}
$$
Recall that
$$
\sigma^2 = \norm{ \sum\nolimits_k \mtx{A}_k^2 } = \lambda_{\max}\left(\sum\nolimits_k \mtx{A}_k^2\right).
$$
Corollary~\ref{cor:main-result-indep} delivers
\begin{equation} \label{eqn:gauss-part-1}
\Prob{ \lambda_{\max}\left( \sum\nolimits_k \xi_k \mtx{A}_k \right) \geq t }
	\leq d \cdot \inf_{\theta > 0} \econst^{ - \theta t + g(\theta) \cdot \sigma^2 }
	= d \cdot \econst^{-t^2/2\sigma^2}.
\end{equation}
For the record, the infimum is attained when $\theta = t/\sigma^2$.

To obtain the norm bound~\eqref{eqn:twosided-matrix-gauss}, recall that $\norm{ \mtx{Y} } = \max\{ \lambda_{\max}(\mtx{Y}), - \lambda_{\min}(\mtx{Y}) \}$.  Standard Gaussian variables and Rademacher variables are symmetric, so the inequality~\eqref{eqn:gauss-part-1} implies
$$
\Prob{ - \lambda_{\min}\left( \sum\nolimits_k \xi_k \mtx{A}_k \right) \geq t }
	= \Prob{ \lambda_{\max}\left( \sum\nolimits_k (-\xi_k) \mtx{A}_k \right) \geq t }
	\leq d \cdot \econst^{-t^2/2\sigma^2}.
$$
Apply the union bound to the estimates for $\lambda_{\max}$ and $-\lambda_{\min}$ to complete the proof.
\end{proof}

The result for a series with rectangular matrix coefficients follows immediately when we apply Theorem~\ref{thm:rad-gauss-series} to the s.a.~dilation of the series.

\begin{proof}[Proof of Corollary~\ref{cor:rect-rad-gauss-series}]
Let $\{\xi_k\}$ be a finite sequence of independent standard normal random variables or independent Rademacher random variables.  Consider the sequence $\{ \xi_k \coll{S}( \mtx{B}_k ) \}$ of random s.a.~matrices with dimension $d_1 + d_2$.
The spectral identity~\eqref{eqn:sa-norm} ensures that
$$
\norm{ \sum\nolimits_k \xi_k \mtx{B}_k }
	= \lambda_{\max}\left( \coll{S}\left( \sum\nolimits_k \xi_k \mtx{B}_k \right) \right)
	= \lambda_{\max}\left(\sum\nolimits_k \xi_k \coll{S}(\mtx{B}_k)\right).
$$
Thus, we may invoke Theorem~\ref{thm:rad-gauss-series} to obtain a probability inequality for the norm of the series.
Simply observe that the matrix variance parameter~\eqref{eqn:matrix-gauss-sd} satisfies the relation
$$
\sigma^2 = \norm{ \sum\nolimits_k \coll{S}(\mtx{B}_k)^2 }
	= \norm{ \begin{bmatrix} \sum\nolimits_k \mtx{B}_k \mtx{B}_k^\adj & \mtx{0} \\
	 \mtx{0} & \sum\nolimits_k \mtx{B}_k^\adj\mtx{B}_k \end{bmatrix} }
	= \max\left\{
	\norm{ \sum\nolimits_k \mtx{B}_k \mtx{B}_k^\adj}, \
	\norm{ \sum\nolimits_k \mtx{B}_k^\adj \mtx{B}_k}
	\right\}
$$
on account of the identity~\eqref{eqn:sa-modulus} for the square of the s.a.~dilation.
\end{proof}

\subsection{Application: A Gaussian Matrix with Nonuniform Variances}

It may not be immediately clear why abstract probability inequalities, such as Theorem~\ref{thm:rad-gauss-series} and Corollary~\ref{cor:rect-rad-gauss-series}, deliver information about interesting random matrices that arise in practice.  
Let us describe a simple application that speaks to this concern.

Fix a $d_1 \times d_2$ matrix $\mtx{B}$, and draw a random $d_1 \times d_2$ matrix $\mtx{\Gamma}$ whose entries are independent standard normal variables.
Let $\odot$ denote the componentwise (i.e., Schur or Hadamard) product of matrices.  Construct the random matrix $\mtx{\Gamma} \odot \mtx{B}$, and observe that its $(j, k)$ component is a Gaussian variable with mean zero and variance $\abssq{b_{jk}}$.  We claim that
\begin{equation} \label{eqn:nonunif-gauss}
\Prob{ \norm{ \mtx{\Gamma} \odot \mtx{B} } \geq t }
	\leq (d_1 + d_2) \cdot \econst^{-t^2/2\sigma^2}
	\quad\text{where}\quad
	\sigma^2 = \max\left\{ \max\nolimits_j \normsq{\vct{b}_{j:}}, \
		\max\nolimits_k \normsq{\vct{b}_{:k}} \right\}.
\end{equation}
The symbols $\vct{b}_{j:}$ and $\vct{b}_{:k}$ represent the $j$th row and $k$th column of the matrix $\mtx{B}$.  An immediate consequence of~\eqref{eqn:nonunif-gauss} is that the median of the norm satisfies
\begin{equation} \label{eqn:nonunif-gauss-med}
\mathbb{M}( \norm{ \mtx{\Gamma} \odot \mtx{B} } ) \leq \sigma \sqrt{2\log( 2(d_1 + d_2))}.
\end{equation}


\noindent
There are nonuniform Gaussian matrices where the estimate~\eqref{eqn:nonunif-gauss-med} for the median has the correct order and other examples where the logarithmic factor is parasitic; see~\S\S\ref{sec:gauss-expect}--\ref{sec:dim-factor} below.
The reader may also wish to juxtapose~\eqref{eqn:nonunif-gauss-med} with the work of Seginer~\cite[Thm.~3.1]{Seg00:Expected-Norm} and Lata{\l}a~\cite[Thm.~1]{Lat05:Some-Estimates} although these results are not fully comparable.

To establish~\eqref{eqn:nonunif-gauss}, we first decompose the matrix of interest as a Gaussian series:
$$
\mtx{\Gamma} \odot \mtx{B} = \sum\nolimits_{jk} \gamma_{jk} \cdot b_{jk} {\bf E}_{jk}.
$$
Next, we must determine the variance parameter.  Note that
$$
\sum\nolimits_{jk} (b_{jk} {\bf E}_{jk})(b_{jk} {\bf E}_{jk})^\adj 
	= \sum\nolimits_{j} \left(\sum\nolimits_k \abssq{b_{jk}}\right) {\bf E}_{jj}
	= \diag( \normsq{ \vct{b}_{1:}}, \normsq{\vct{b}_{2:}},
	\dots, \normsq{\vct{b}_{d_1:}} ).
$$
Similarly,
$$
\sum\nolimits_{jk} (b_{jk} {\bf E}_{jk})^\adj(b_{jk} {\bf E}_{jk})
	= \sum\nolimits_{k} \left(\sum\nolimits_j \abssq{b_{jk}}\right) {\bf E}_{kk}
	= \diag( \normsq{ \vct{b}_{:1}}, \normsq{\vct{b}_{:2}},
	\dots, \normsq{\vct{b}_{:d_2}} ).
$$
Therefore,
\begin{align*}
\sigma^2 &= \max\left\{ \smnorm{}{ \diag( \normsq{ \vct{b}_{1:}}, \normsq{\vct{b}_{2:}},
	\dots, \normsq{\vct{b}_{d_1:}} ) }, \
	\smnorm{}{ \diag( \normsq{ \vct{b}_{:1}}, \normsq{\vct{b}_{:2}},
	\dots, \normsq{\vct{b}_{:d_2}} ) } \right\} \\
	&= \max\left\{ \max\nolimits_j \normsq{\vct{b}_{j:}}, \
		\max\nolimits_k \normsq{\vct{b}_{:k}} \right\}.
\end{align*}
An application of Corollary~\ref{cor:rect-rad-gauss-series} yields the tail bound~\eqref{eqn:nonunif-gauss}.  



\subsection{Controlling the Expectation} \label{sec:gauss-expect}

A remarkable feature of Theorem~\ref{thm:rad-gauss-series} is that it always allows us to obtain reasonably accurate estimates for the expected norm of the s.a.~Gaussian series
\begin{equation} \label{eqn:gauss-series}
\mtx{Y} = \sum\nolimits_{k} \gamma_k \mtx{A}_k.
\end{equation}
To establish this point, we first compute upper and lower bounds for the second moment of $\norm{\mtx{Y}}$.  Theorem~\ref{thm:rad-gauss-series} yields
\begin{equation*}
	\Expect( \normsq{ \mtx{Y} } )
	= \int_0^\infty \Prob{ \norm{ \mtx{Y} } > \sqrt{t} } \idiff{t}
	\leq 2 \sigma^2 \log(2d) + 2d \int_{2 \sigma^2 \log(2d)}^\infty
	\econst^{-t/2\sigma^2} \idiff{t}
	= 2\sigma^2 \log(2 \econst d).
\end{equation*}
Jensen's inequality furnishes the lower estimate:
$$
\Expect( \normsq{\mtx{Y}} )
	= \Expect \norm{ \mtx{Y}^2 }
	\geq \norm{ \Expect (\mtx{Y}^2) }
	= \norm{ \sum\nolimits_k \mtx{A}_k^2 }
	= \sigma^2.
$$
The (homogeneous) first and second moment of the norm of a Gaussian series 
are equivalent up to a universal constant~\cite[Cor.~3.2]{LT91:Probability-Banach},
so we conclude that
\begin{equation} \label{eqn:sigma-mean}
\cnst{c} \sigma \leq \Expect \norm{ \mtx{Y} } \leq \sigma \sqrt{2 \log(2\econst d)}.
\end{equation}
This argument demonstrates that the matrix variance parameter $\sigma^2$ controls the expected norm $\Expect \norm{ \mtx{Y} }$ up to a factor that depends very weakly on the dimension.  A similar remark applies to the median value $\mathbb{M}( \norm{\mtx{Y}} )$.

\subsection{The Dimensional Factor} \label{sec:dim-factor}


In the inequality~\eqref{eqn:sigma-mean}, the gap between the upper and lower bounds for $\Expect \norm{\mtx{Y}}$ arises because of the dimensional factor $d$ in the statement~\eqref{eqn:twosided-matrix-gauss}.  This dimensional dependence is a new feature of probability inequalities in the matrix setting.  The extra term appears in each of our main results, and it is usually possible to identify a simple case where it is necessary.


In particular, we cannot remove the factor $d$ from the probability bound in Theorem~\ref{thm:rad-gauss-series}.  Observe that the norm of a diagonal Gaussian matrix is typically bounded below:
$$
\norm{ \sum\nolimits_{k=1}^d \gamma_k {\bf E}_{kk} }
	= \max\nolimits_k \abs{\gamma_k}
	> \sqrt{2\log d}
\quad\text{with high probability.}	
$$
Theorem~\ref{thm:rad-gauss-series} delivers the following tail bound for this series.
$$
\Prob{ \norm{ \sum\nolimits_{k=1}^d \gamma_k {\bf E}_{kk} }
	\geq t } \leq 2d \cdot \econst^{-t^2/2}.
$$
The factor $2d$ ensures that this probability inequality does not become effective until $t \geq \sqrt{2 \log(2d)}$, \lang{comme il faut}.

We can also identify situations where the dimensional term produces an overestimate of the expected norm.  For instance, consider a $d$-dimensional matrix drawn from the unnormalized Gaussian orthogonal ensemble (GOE): 
$$
\mtx{W} = \sum_{1 \leq j \leq k \leq d} \gamma_{jk} ({\bf E}_{jk} + {\bf E}_{kj})
$$
The literature contains a sharp bound for the expected norm of this matrix:
\begin{equation} \label{eqn:goe-bd}
\Expect \norm{ \mtx{W} } \leq 2\sqrt{d}
\end{equation}
The result~\eqref{eqn:goe-bd} follows from ideas of Gordon~\cite{Gor85:Some-Inequalities,Gor92:Majorization-Gaussian} elaborated in~\cite[Thm.~2.11]{DS02:Local-Operator}. 
Meanwhile, integrating the tail bound~\eqref{eqn:twosided-matrix-gauss} from Theorem~\ref{thm:rad-gauss-series} yields the weaker result
\begin{equation} \label{eqn:lieb-goe-bd}
\Expect \norm{ \mtx{W} } \leq \sqrt{(d + 3) \log(2\econst d)}.
\end{equation}
The estimate~\eqref{eqn:lieb-goe-bd} is too large by a factor of about $\sqrt{\log d}$, which is the worst possible discrepancy in view of~\eqref{eqn:sigma-mean}.

\begin{rem}[Effective Dimension]
Let us stress that the \emph{nominal} dimension of the matrices does not play a role in Theorem~\ref{thm:rad-gauss-series}.  If the ranges of the matrices $\mtx{A}_1, \mtx{A}_2, \dots$ are contained within a fixed $r$-dimensional subspace, we can replace the ambient dimension~$d$ with the effective dimension~$r$.  A similar remark applies to our other results.
\end{rem}



\subsection{Comparison with Concentration Inequalities}


It is fruitful to think about Theorem~\ref{thm:rad-gauss-series} as a statement that the matrix Gaussian series~\eqref{eqn:gauss-series} typically falls near its expectation \emph{as a random matrix} when we measure the size of deviations using the operator norm:
\begin{equation} \label{eqn:matrix-gauss-concentration}
\Prob{ \norm{ \mtx{Y} - \Expect \mtx{Y} } \geq t } \leq 2d \cdot \econst^{-t^2/2\sigma^2}.
\end{equation}
In contrast, the classical concentration inequality~\cite[Thm.~1.7.6]{Bog98:Gaussian-Measures} concerns the variation of \emph{the norm} about its mean value:
\begin{equation} \label{eqn:gauss-concentration}
\Prob{ \abs{ \, \norm{ \mtx{Y} }
	- \Expect \norm{ \mtx{Y} } \, } \geq t }
	\leq 2 \cdot \econst^{- t^2 / 2 \sigma_{*}^2}
\end{equation}
where the scale for deviations depends on the \emph{weak variance} parameter
\begin{equation} \label{eqn:gauss-weak-var}
\sigma_{*}^2 := \sup \left\{ 
\sum\nolimits_k \abssq{ \vct{u}^\adj \mtx{A}_k \vct{v} } :
\norm{\vct{u}} = \norm{\vct{v}} = 1 \right\}.
\end{equation}
It can be shown~\cite[Cor.~3.2]{LT91:Probability-Banach} that the bound~\eqref{eqn:gauss-concentration} is asymptotically sharp as $t \to \infty$.



Let us elaborate on the relationship between the matrix variance $\sigma^2$ defined in~\eqref{eqn:matrix-gauss-sd} and the weak variance $\sigma_{*}^2$ appearing in~\eqref{eqn:gauss-weak-var}.  First, note that
\begin{equation} \label{eqn:weak-variance-ub}
\sigma_{*}^2
	\leq \sup_{\norm{\vct{u}} = 1} \sum\nolimits_k \vct{u}^\adj \mtx{A}_k^2 \, \vct{u}
	= \norm{ \sum\nolimits_k \mtx{A}_k^2 }
	= \sigma^2.
\end{equation}
Equality holds in~\eqref{eqn:weak-variance-ub} when, for example, the family $\{\mtx{A}_k\}$ commutes.  We can also establish a reverse inequality.
\begin{equation} \label{eqn:weak-variance-lb}
\sigma^2 = \norm{ \sum\nolimits_k \mtx{A}_k \left( \sum\nolimits_j \onevct_j \onevct_j^\adj \right) \mtx{A}_k }
	\leq \sum\nolimits_j \sup_{\norm{\vct{u}} = 1} \sum\nolimits_k
	\abssq{ \vct{u}^\adj \mtx{A}_k \onevct_j }
	\leq d \cdot \sigma_*^2
\end{equation}
where $\{\onevct_j : j = 1, \dots, d \}$ is the standard basis for $\Rspace{d}$.  In the worst case%
\footnote{A worst-case example occurs with high probability when the sequence $\{ \mtx{A}_k : k = 1, \dots, d \}$ consists of independent matrices drawn from the $d$-dimensional GOE, but the proof seems to be complicated.},
the bound~\eqref{eqn:weak-variance-lb} has roughly the correct order.



In summary, the matrix concentration inequality~\eqref{eqn:matrix-gauss-concentration} always leads to a good estimate for the expected norm $\Expect \norm{ \mtx{Y} }$.  Nevertheless, the presence of the parameter $\sigma^2$ in the tail bound can lead to a significant overestimate of the probability that $\norm{\mtx{Y}}$ is large.  On the other hand, the classical inequality~\eqref{eqn:gauss-concentration} contains no information about the mean, but it always produces a sharp large-deviation bound.  Therefore, the two results complement each other well.

\subsection{Noncommutative Moment Inequalities} \label{sec:nc-moments}

The matrix Laplace transform bound, Proposition~\ref{prop:laplace-transform} demonstrates that we can bound tail probabilities for the norm of a random series by controlling the matrix mgf.
In certain special cases, it is possible to bound the matrix mgf using noncommutative (nc) moment inequalities.  Let us describe how to establish Theorem~\ref{thm:rad-gauss-series} in this fashion.  This material is unrelated to the main development, so the reader may skip it with impunity.  

The nc Khintchine inequality provides an estimate for the expectation of the $(2p)$th moment of the Schatten $2p$-norm of a matrix Gaussian series~\cite{L-P86:Inegalites-Khintchine,LPP91:Noncommutative-Khintchine,Pis03:Introduction-Operator}.  The most elementary formulation of this result states that
\begin{equation} \label{eqn:nc-Khintchine}
\Expect \trace \left( \sum\nolimits_k \gamma_k \mtx{A}_k \right)^{2p}
	\leq \cnst{C}_{2p} \cdot \trace \left( \sum\nolimits_k \mtx{A}_k^2 \right)^p
	\quad\text{for $p = 1, 2, 3, \dots$.}
\end{equation}
Buchholz~\cite[Thm.~5]{Buc01:Operator-Khintchine} has shown that the optimal constant in~\eqref{eqn:nc-Khintchine} satisfies
$$
\cnst{C}_{2p} := \Expect \abs{\gamma_1}^{2p}
	= (2p-1)!!
	= \frac{(2p)!}{p! \, 2^p}.
$$
The bound~\eqref{eqn:nc-Khintchine} also holds with the same constant when we replace $\{\gamma_k\}$ by a sequence of independent Rademacher variables~\cite[Thm.~5]{Buc05:Optimal-Constants}.  


The family~\eqref{eqn:nc-Khintchine} of inequalities allows us to develop a short proof of the tail bound for matrix Gaussian and Rademacher series.

\begin{proof}[Alternative Proof of Theorem~\ref{thm:rad-gauss-series}]
Proposition~\ref{prop:laplace-transform} yields
\begin{equation} \label{eqn:gauss-lt}
\Prob{ \lambda_{\max}\left( \sum\nolimits_k \gamma_k \mtx{A}_k \right) \geq t }
	\leq \inf_{\theta > 0} \left\{ \econst^{-\theta t}
	\cdot \Expect \trace \exp\left( \theta \sum\nolimits_k \gamma_k \mtx{A}_k \right) \right\}.
\end{equation}
We may use~\eqref{eqn:nc-Khintchine} to bound the Taylor series for the matrix mgf term by term:
\begin{align} \label{eqn:gauss-moment-bd}
\Expect \trace \exp\left( \theta \sum\nolimits_k \gamma_k \mtx{A}_k \right)	
	&= \sum\nolimits_{p = 0}^\infty \frac{\theta^{2p}}{(2p)!}
		\Expect \trace \left( \sum\nolimits_k \gamma_k \mtx{A}_k \right)^{2p} \notag \\
	&\leq \sum\nolimits_{p = 0}^\infty \frac{\theta^{2p}}{p! \, 2^p}
		\trace \left( \sum\nolimits_k \mtx{A}_k^2 \right)^{p}
	= \trace \exp\left( \frac{\theta^2}{2} \sum\nolimits_k \mtx{A}_k^2 \right).
\end{align}
Substitute~\eqref{eqn:gauss-moment-bd} into~\eqref{eqn:gauss-lt}, and select $\theta = t/\sigma^2$ to complete the minimization.
\end{proof}

We may regard the mgf bound~\eqref{eqn:gauss-moment-bd} as an ``exponential generating function'' for the family of nc Khintchine inequalities~\eqref{eqn:nc-Khintchine}, but---unfortunately---the nc Khintchine inequalities \emph{do not} follow as a consequence of this mgf bound. 
Recall that Lieb's result, Theorem~\ref{thm:lieb}, also delivers a proof of the inequality~\eqref{eqn:gauss-moment-bd}.  This observation suggests that it might be possible to use Lieb's theorem to prove the nc Khintchine inequalities~\eqref{eqn:nc-Khintchine}.  We regard this as a tantalizing open question.

\subsection{Comparison with the Ahlswede--Winter Bound} \label{sec:aw-gauss}

In~\S\ref{sec:aw}, we describe how Ahlswede and Winter go about bounding the matrix mgf~\cite[App.]{AW02:Strong-Converse}.  It is natural to ask how inequalities developed using their approach compare with the results in this paper.

Gaussian series provide an excellent illustration of the discrepancy between the two techniques.  In this case, the Ahlswede--Winter method yields the probability inequality
\begin{equation} \label{eqn:aw-matrix-gauss}
\Prob{ \norm{ \sum\nolimits_k \gamma_k \mtx{A}_k } \geq t }
	\leq 2d \cdot \econst^{- t^2 / 2\sigma_{\rm AW}^2 }
\quad\text{where $\sigma_{\rm AW}^2 := \sum\nolimits_k \norm{\mtx{A}_k^2}$.}
\end{equation}
The estimate~\eqref{eqn:aw-matrix-gauss} should be compared with our bound~\eqref{eqn:twosided-matrix-gauss}.
The Ahlswede--Winter variance parameter $\sigma_{\rm AW}^2$
always dominates the matrix variance parameter~\eqref{eqn:matrix-gauss-sd} because
$$
\sigma^2 = \norm{ \sum\nolimits_k \mtx{A}_k^2 }
	\leq \sum\nolimits_k \norm{\mtx{A}_k^2}
	= \sigma^2_{\rm AW}.
$$
The two variance parameters rarely coincide, and the best reverse inequality is
$$
\sigma_{\rm AW}^2 \leq \sum\nolimits_k \trace \mtx{A}_k^2
	\leq d \cdot \norm{ \sum\nolimits_k \mtx{A}_k^2 }
	= d \cdot \sigma^2.
$$
This 
worst-case behavior is typical. For instance, consider the two Gaussian matrices presented in~\S\ref{sec:dim-factor}.  The Ahlswede--Winter tail bound~\eqref{eqn:aw-matrix-gauss} provides essentially no information about the norm of either matrix.


\begin{rem}[Moment Inequalities]
There is an alternative approach to establishing the result~\eqref{eqn:aw-matrix-gauss} that parallels the method presented in~\S\ref{sec:nc-moments}.  We simply bound the Taylor series of the matrix mgf term by term using an appropriate family of moment inequalities:
$$
\Expect \trace \left( \sum\nolimits_k \gamma_k \mtx{A}_k \right)^{2p}
	\leq \cnst{C}_{2p} \cdot \left( \sum\nolimits_k \left[ \trace(\mtx{A}_k^{2p}) \right]^{1/p} \right)^p
\quad\text{where}\quad
\cnst{C}_{2p} := \frac{(2p)!}{p! \, 2^p}
\quad\text{for $p = 1, 2, 3, \dots$.}
$$
These estimates follow from a result of Tomczak--Jaegermann~\cite[Thm.~3.1]{TJ74:Moduli-Smoothness} for Rademacher series together with the central limit theorem.
\end{rem}

\section{Sums of Random Positive-Semidefinite Matrices} \label{sec:chernoff}

The classical Chernoff bounds concern the sum of independent, nonnegative, and uniformly bounded random variables.  In sympathy, matrix Chernoff bounds describe the extreme eigenvalues of a sum of independent, psd random matrices whose maximum eigenvalues are subject to a uniform bound.  These probability inequalities demonstrate that the upper and lower tails of the sum exhibit binomial-type behavior.

Our first result parallels the strongest versions of the scalar Chernoff inequality for the proportion of successes in a sequence of independent (but not identical) Bernoulli trials~\cite[Exer.~7]{Lug09:Concentration-Measure}. 

\begin{thm}[Matrix Chernoff I] \label{thm:chernoff-full}
Consider a sequence $\{ \mtx{X}_k : k = 1, 2, \dots, n \}$ of independent, random, self-adjoint matrices that satisfy
$$
\mtx{X}_k \psdge \mtx{0}
\quad\text{and}\quad
\lambda_{\max}(\mtx{X}_k) \leq 1
\quad\text{almost surely}.
$$
Compute the minimum and maximum eigenvalues of the average expectation,
$$
\bar{\mu}_{\min} := \lambda_{\min}\left( \frac{1}{n}  \sum\nolimits_{k=1}^n \Expect \mtx{X}_k \right)
\quad\text{and}\quad
\bar{\mu}_{\max} := \lambda_{\max}\left( \frac{1}{n}  \sum\nolimits_{k=1}^n \Expect \mtx{X}_k \right).
$$
Then
\begin{align*}
\Prob{ \lambda_{\min}\left(\frac{1}{n} \sum\nolimits_{k=1}^n \mtx{X}_k \right) \leq \alpha }
	&\leq d \cdot \econst^{-n \cdot {\rm D}( \alpha \, \Vert \, \bar{\mu}_{\min} )} 
	\quad\text{for $0 \leq \alpha \leq \bar{\mu}_{\min}$, and} \\ 
\Prob{ \lambda_{\max}\left(\frac{1}{n} \sum\nolimits_{k=1}^n \mtx{X}_k\right) \geq \alpha }
	&\leq d \cdot \econst^{-n \cdot {\rm D}( \alpha \, \Vert \, \bar{\mu}_{\max} )}
	\quad\text{for $\bar{\mu}_{\max} \leq \alpha \leq 1$}. 
\end{align*}
The binary information divergence 
${\rm D}( a \, \Vert \, u ) := a(\log(a) - \log(u)) + (1-a)(\log(1-a) - \log(1-u))$
for $a,u \in [0,1]$.
\end{thm}


We have found that the following weaker version of Theorem~\ref{thm:chernoff-full} produces excellent results but is simpler to apply.  This corollary corresponds with the usual statement of the scalar Chernoff inequalities for sums of nonnegative random variables; see~\cite[Exer.~8]{Lug09:Concentration-Measure} or~\cite[\S4.1]{MR95:Randomized-Algorithms}.

\begin{cor}[Matrix Chernoff II] \label{cor:chernoff}
Consider a finite sequence $\{ \mtx{X}_k \}$ of independent, random, self-adjoint matrices that satisfy
$$
\mtx{X}_k \psdge \mtx{0}
\quad\text{and}\quad
\lambda_{\max}(\mtx{X}_k) \leq R
\quad\text{almost surely}.
$$
Compute the minimum and maximum eigenvalues of the sum of expectations,
$$
\mu_{\min} := \lambda_{\min}\left( \sum\nolimits_{k} \Expect \mtx{X}_k \right)
\quad\text{and}\quad
\mu_{\max} := \lambda_{\max}\left( \sum\nolimits_{k} \Expect \mtx{X}_k \right).
$$
Then
\begin{align*}
\Prob{ \lambda_{\min}\left(\sum\nolimits_k \mtx{X}_k \right) \leq (1 - \delta) \mu_{\min} }
	&\leq d \cdot \left[ \frac{\econst^{-\delta}}{(1 - \delta)^{1-\delta}} \right]^{\mu_{\min}/R} 
	\quad\text{for $\delta \in [0, 1]$, and} \\
\Prob{ \lambda_{\max}\left( \sum\nolimits_k \mtx{X}_k \right) \geq (1 + \delta) \mu_{\max} }
	&\leq d \cdot \left[ \frac{\econst^{\delta}}{(1 + \delta)^{1+\delta}} \right]^{\mu_{\max} / R} 
	\quad\text{for $\delta \geq 0$.} 
\end{align*}
\end{cor}

The proofs of Theorem~\ref{thm:chernoff-full} and Corollary~\ref{cor:chernoff} appear below in Section~\ref{sec:chernoff-proofs}.  We continue this discussion with some telegraphic remarks concerning various aspects of the Chernoff bounds.



\begin{rem}[Related Inequalities]
The following standard simplification of Corollary~\ref{cor:chernoff} is useful. 
\begin{align*}
\Prob{ \lambda_{\min}\left(\sum\nolimits_k \mtx{X}_k \right)
	\leq t \mu_{\min} }
	&\leq d \cdot \econst^{-(1-t)^2 \mu_{\min}/2R}
	\quad\text{for $t \in [0, 1]$, and} \\
\Prob{ \lambda_{\max}\left(\sum\nolimits_{k} \mtx{X}_k \right)
	\geq t \mu_{\max} }
	&\leq d \cdot \left[ \frac{\econst}{t} \right]^{t \mu_{\max}/R}
\quad\text{for $t \geq \econst$.}
\end{align*}
These inequalities manifest that the minimum eigenvalue has normal-type behavior and the maximum eigenvalue exhibits Poisson-type decay.

\end{rem}

\begin{rem}[Applications]
Matrix Chernoff inequalities are very effective for studying random matrices with independent columns.  Consider a rectangular random matrix
$$
\mtx{Z} = \begin{bmatrix} \vct{z}_1 & \vct{z}_2 & \dots & \vct{z}_n \end{bmatrix}
$$
where $\{ \vct{z}_k \}$ is a family of independent random vectors in $\mathbb{C}^m$.
The norm of $\mtx{Z}$ satisfies
$$
\norm{ \mtx{Z} }^2 = \lambda_{\max}( \mtx{ZZ}^\adj )
	= \lambda_{\max}\left( \sum\nolimits_{k=1}^n \vct{z}_k \vct{z}_k^\adj \right).
$$
Similarly, the minimum singular value $s_m$ of the matrix satisfies
$$
s_{m}( \mtx{Z} )^2 = \lambda_{\min}( \mtx{ZZ}^\adj )
	= \lambda_{\min}\left( \sum\nolimits_{k=1}^n \vct{z}_k \vct{z}_k^\adj \right).
$$
In each case, the summands are stochastically independent and psd, so the matrix Chernoff bounds apply.  See~\cite{Tro10:Improved-Analysis} for a problem where this method applies.
\end{rem}

\begin{rem}[Expectations]
Corollary~\ref{cor:chernoff} produces accurate estimates for the expectation of the maximum eigenvalue:
$$
\mu_{\max} \leq \Expect \lambda_{\max}\left(\sum\nolimits_k \mtx{X}_k \right)
	\leq \cnst{C} \cdot \max\left\{ \mu_{\max},\ R \log d \right\}.
$$
The lower bound is Jensen's inequality; the upper bound follows from a messy---but standard---calculation.  Observe that the dimensional dependence vanishes when the mean $\mu_{\max}$ is sufficiently large in comparison with the upper bound $R$!
\end{rem}

\begin{rem}[Dimensional Factor]
The factor $d$ in the Chernoff bounds cannot be omitted because of the coupon collector's problem~\cite[\S3.6]{MR95:Randomized-Algorithms}.  Consider a $d$-dimensional random matrix $\mtx{X}$ with the distribution
$$
\mtx{X} = {\bf E}_{jj} \quad\text{with probability $d^{-1}$ for each $j = 1, 2, \dots d$.}
$$
If $\{\mtx{X}_k\}$ is a sequence of independent random matrices with the same distribution as $\mtx{X}$, then
$$
\lambda_{\min}\left(\sum\nolimits_{k=1}^n \mtx{X}_k \right)
	= 0
\quad\text{with high probability unless $n > d\log d$.}
$$
The dimensional factor in the lower Chernoff bound reflects this fact.  The same example shows that the upper Chernoff bound must also exhibit a dimensional dependence.  We have extracted this idea from~\cite[Sec.~3.5]{RV07:Sampling-Large}.
\end{rem}


\begin{rem}[Previous Work]
Theorem~\ref{thm:chernoff-full} is a considerable strengthening of the matrix Chernoff bound established by Ahlswede and Winter~\cite[Thm.~19]{AW02:Strong-Converse}.  Their proof requires the extra assumption that the summands are identically distributed, in which case their result matches Theorem~\ref{thm:chernoff-full}.
\end{rem}

\subsection{Proofs} \label{sec:chernoff-proofs}

To establish the matrix Chernoff inequalities, we commence with a semidefinite bound for the matrix mgf of a random psd contraction.



\begin{lemma}[Chernoff mgf] \label{lem:chernoff-mgf}
Suppose that $\mtx{X}$ is a random psd~matrix that satisfies $\lambda_{\max}(\mtx{X}) \leq 1$.  Then
$$
\Expect \econst^{\theta\mtx{X}}
	\psdle \Id +  (\econst^{\theta} - 1)(\Expect \mtx{X})
\quad\text{for $\theta \in \mathbb{R}$}.
$$
\end{lemma}

The proof of Lemma~\ref{lem:chernoff-mgf} parallels the classical argument; the matrix adaptation is due to Ahlswede and Winter~\cite[Thm.~19]{AW02:Strong-Converse}.

\begin{proof}
Consider the function $f(x) = \econst^{\theta x}$.  Since $f$ is convex, its graph lies below the chord connecting two points.  In particular,
$$
f(x) \leq f(0) + [f(1) - f(0)] \cdot x
\quad\text{for $x \in [0, 1]$.}
$$
More explicitly,
$$
\econst^{\theta x} \leq 1 + (\econst^{\theta} - 1) \cdot x
\quad\text{for $x \in [0, 1]$.}
$$
The eigenvalues of $\mtx{X}$ lie in the interval $[0, 1]$, so the transfer rule~\eqref{eqn:scalar-matrix} implies that
$$
\econst^{\theta \mtx{X}} \psdle \Id + (\econst^{\theta} - 1) \mtx{X}.
$$
Expectation respects the semidefinite order, so
$$
\Expect \econst^{\theta \mtx{X}}
	\psdle \Id + (\econst^{\theta} - 1) (\Expect\mtx{X}).
$$
This is the advertised conclusion.
\end{proof}

We prove the upper Chernoff bounds first because the argument is slightly easier.

\begin{proof}[Proof of Theorem~\ref{thm:chernoff-full}, Upper Bound]
The Chernoff mgf bound, Lemma~\ref{lem:chernoff-mgf}, states that
$$
\Expect \econst^{\theta\mtx{X}_k}
	\psdle \Id +  g(\theta) \cdot (\Expect \mtx{X}_k)
	\quad\text{where $g(\theta) := \econst^{\theta} - 1$ for $\theta > 0$.}
$$
As a result, Corollary~\ref{cor:mgf-master-ineq} implies
\begin{align} \label{eqn:upper-chernoff-pf}
\Prob{ \lambda_{\max}\left(\sum\nolimits_k \mtx{X}_k \right) \geq t }
	&\leq d \cdot \exp\left( -\theta t +
	n \cdot \log \lambda_{\max}\left( \frac{1}{n} \sum\nolimits_k (\Id + g(\theta) \cdot
	\Expect \mtx{X}_k) \right) \right) \notag \\
	&= d \cdot \exp\left( -\theta t +
	n \cdot \log \lambda_{\max}\left( \Id + g(\theta) \cdot
	\frac{1}{n} \sum\nolimits_k \Expect \mtx{X}_k \right) \right) \notag \\
	&= d \cdot \exp\left( -\theta t +
	n \cdot \log\left( 1 + g(\theta) \cdot \bar{\mu}_{\max} \right) \right).
\end{align}
The third relation follows from basic properties of the eigenvalue map and the definition of $\bar{\mu}_{\max}$.
Make the change of variables $t \mapsto n \alpha$.  The right-hand side is smallest when
$$
\theta = \log(\alpha / (1 - \alpha)) - \log(\bar{\mu}_{\max}/(1-\bar{\mu}_{\max})).
$$
Substitute these quantities into~\eqref{eqn:upper-chernoff-pf} to obtain the information divergence upper bound.
\end{proof}

\begin{proof}[Proof of Corollary~\ref{cor:chernoff}, Upper Bound]
Assume that the summands satisfy the uniform eigenvalue bound with $R = 1$; the general result follows by re-scaling.  The shortest route to the weaker Chernoff upper bound starts at~\eqref{eqn:upper-chernoff-pf}.  The numerical inequality $\log(1 + x) \leq x$, valid for $x > -1$, implies that 
$$
\Prob{ \lambda_{\max}\left(\sum\nolimits_k \mtx{X}_k \right) \geq t}
	\leq d \cdot \exp\left( -\theta t +
	g(\theta) \cdot n \bar{\mu}_{\max} \right)
	= d \cdot \exp\left( -\theta t +
	g(\theta) \cdot \mu_{\max} \right)
$$
Make the change of variables $t \mapsto (1+\delta) \mu_{\max}$, and select the parameter $\theta = \log(1 + \delta)$.  Simplify the resulting tail bound to complete the proof. 
\end{proof}

The lower bounds follow from a closely related argument.

\begin{proof}[Proof of Theorem~\ref{thm:chernoff-full}, Lower Bound]
We intend to apply Corollary~\ref{cor:mgf-master-ineq} to the sequence $\{-\mtx{X}_k\}$.  In this case, the Chernoff mgf, Lemma~\ref{lem:chernoff-mgf}, states that
$$
\Expect \econst^{\theta(-\mtx{X}_k)}
	= \Expect \econst^{(-\theta)\mtx{X}_k}
	\psdle \Id - g(\theta) \cdot (\Expect \mtx{X}_k)
	\quad\text{where $g(\theta) := 1 - \econst^{-\theta}$ for $\theta > 0$.}
$$
The minimum eigenvalue $\lambda_{\min}(-\mtx{A}) = - \lambda_{\max}(\mtx{A})$, so we can apply Corollary~\ref{cor:mgf-master-ineq} as follows.
\begin{align} \label{eqn:lower-chernoff-pf}
\Prob{ \lambda_{\min}\left(\sum\nolimits_k \mtx{X}_k \right) \leq t }
	&= \Prob{ \lambda_{\max}\left(\sum\nolimits_k (- \mtx{X}_k) \right) \geq -t } \notag \\
	&\leq d \cdot \exp\left( \theta t +
 	n \cdot \log \lambda_{\max}\left( \frac{1}{n} \sum\nolimits_k ( \Id - g(\theta) \cdot
	\Expect \mtx{X}_k) \right) \right) \notag \\
	&= d \cdot \exp\left( \theta t +
	n \cdot \log\left( 1 - g(\theta) \cdot \lambda_{\min}\left(
	\frac{1}{n} \sum\nolimits_{k=1}^n \Expect \mtx{X}_k \right) \right) \right) \notag \\
	&= d \cdot \exp\left( \theta t + 
	n \cdot \log\left( 1 - g(\theta) \cdot \bar{\mu}_{\min} \right) \right).
\end{align}
Make the substitution $t \mapsto n \alpha$. The right-hand side is minimal when
$$
\theta = \log(\bar{\mu}_{\min}/(1-\bar{\mu}_{\min})) - \log(\alpha / (1 - \alpha)).
$$
These steps result in the information divergence lower bound.
\end{proof}

\begin{proof}[Proof of Corollary~\ref{cor:chernoff}, Lower Bound]
As before, assume that the uniform bound $R = 1$.  We obtain the weaker lower bound as a consequence of~\eqref{eqn:lower-chernoff-pf}.  The inequality $\log(1+x) \leq x$ holds for $x > -1$, so we have
$$
\Prob{ \lambda_{\min}\left(\sum\nolimits_k \mtx{X}_k \right) \leq t }
	\leq d \cdot \exp\left( \theta t - g(\theta) \cdot n \bar{\mu}_{\min} \right)
	= d \cdot \exp\left( \theta t - g(\theta) \cdot \mu_{\min} \right)
$$
Make the replacement $t \mapsto (1 -\delta) \mu_{\min}$, and select $\theta = -\log(1-\delta)$ to complete the proof. 
\end{proof}

\begin{rem}[Alternative Proof]
Corollary~\ref{cor:chernoff} can also be established directly using Corollary~\ref{cor:main-result-indep} instead of Corollary~\ref{cor:mgf-master-ineq}.  In this case, we use the mgf bound
$$
\Expect \econst^{\theta \mtx{X}}
	\psdle \exp\left( (\econst^{\theta} - 1) (\Expect \mtx{X}) \right)
	\quad\text{for $\theta \in \mathbb{R}$},
$$
which follows instantly from Lemma~\ref{lem:chernoff-mgf} and the semidefinite relation~\eqref{eqn:Id+A}.  The remaining details mirror the arguments here.
\end{rem}

\section{Matrix Bennett and Bernstein Inequalities} \label{sec:bennett}


In the scalar setting, Bennett and Bernstein inequalities describe the upper tail of a sum of independent, zero-mean random variables that are either bounded or subexponential.  In the matrix case, the analogous results concern a sum of zero-mean random matrices.



Our first result describes the case where the maximum eigenvalue of each summand satisfies a uniform bound.

\begin{thm}[Matrix Bernstein: Bounded Case] \label{thm:bernstein-bdd}
Consider a finite sequence $\{ \mtx{X}_k \}$ of independent, random, self-adjoint matrices with dimension $d$.  Assume that
$$
\Expect \mtx{X}_k = \mtx{0}
\quad\text{and}\quad
\lambda_{\max}(\mtx{X}_k) \leq R
\quad\text{almost surely.}
$$
Compute the norm of the total variance,
$$
\sigma^2 := \norm{ \sum\nolimits_k \Expect \left(\mtx{X}_k^2 \right) }.
$$
Then the following chain of inequalities holds for all $t \geq 0$.
\begin{align*}
\Prob{ \lambda_{\max}\left( \sum\nolimits_k \mtx{X}_k \right) \geq t }
	&\leq d \cdot \exp\left( - \frac{\sigma^2}{R^2} \cdot h\!\left( \frac{Rt}{\sigma^2} \right) \right)
	\tag{i} \\
	&\leq d \cdot \exp\left( \frac{-t^2/2}{\sigma^2 + Rt/3} \right)
	\tag{ii} \\
	&\leq \begin{cases} d \cdot \exp( - 3t^2/8\sigma^2 ) &
	\quad\text{for $t \leq \sigma^2 / R$;} \\
	d \cdot \exp( -3t/8R ) &
	\quad\text{for $t \geq \sigma^2 / R$}.
	\end{cases} \tag{iii}
\end{align*}
The function $h(u) := (1+u)\log(1+u) - u$ for $u \geq 0$.
\end{thm}

Observe that Theorem~\ref{thm:bernstein-bdd} places no assumption on the minimum eigenvalues of the summands, which may be arbitrarily small.  As a consequence, when we apply the result to the two sequences $\{ \mtx{X}_k \}$ and $\{ - \mtx{X}_k \}$, the parameter $R$ may differ.

Theorem~\ref{thm:bernstein-bdd}(i) can be viewed as a matrix version of the Bennett inequality~\cite[Thm.~5]{Lug09:Concentration-Measure}, which implies that the tail probabilities exhibit Poisson-type decay.  Part (ii) parallels a well-known result~\cite[Thm.~6]{Lug09:Concentration-Measure}, which is perhaps the most famous among the probability inequalities attributed to Bernstein.  Part (iii), which we call the split Bernstein inequality, clearly delineates between the normal behavior that occurs at moderate deviations and the slower decay that emerges in the tail.

%

A related inequality holds when we allow the moments of the random matrices to grow at a limited rate, which we interpret as a matrix extension of the moment behavior of a subexponential random variable~\cite[Lem.~4.1.9]{PG02:Decoupling}.

\begin{thm}[Matrix Bernstein: Subexponential Case] \label{thm:bernstein-subexp}
Consider a finite sequence $\{ \mtx{X}_k \}$ of independent, random, self-adjoint matrices with dimension $d$.  Assume that
$$
\Expect \mtx{X}_k = \mtx{0}
\quad\text{and}\quad
\Expect( \mtx{X}_k^p ) \psdle \frac{p!}{2} \cdot R^{p-2} \mtx{A}_k^2
\quad\text{for $p = 2, 3, 4, \dots$.}
$$
Compute the variance parameter
$$
\sigma^2 := \norm{ \sum\nolimits_k \mtx{A}_k^2 }.
$$
Then the following chain of inequalities holds for all $t \geq 0$.
\begin{align*}
\Prob{ \lambda_{\max}\left( \sum\nolimits_k \mtx{X}_k \right) \geq t }
	&\leq d \cdot \exp \left( \frac{-t^2/2}{\sigma^2 + Rt} \right)
		\tag{i} \\
	&\leq \begin{cases}
		d \cdot \exp( -t^2/4\sigma^2 ) &
		\text{for $t \leq \sigma^2 / R$;} \\
		d \cdot \exp( -t/4R  ) &
		\text{for $t \geq \sigma^2 / R$.}
	\end{cases} \tag{ii}
\end{align*}
\end{thm}

The hypotheses of Theorem~\ref{thm:bernstein-subexp} are not fully comparable with the hypotheses of Theorem~\ref{thm:bernstein-bdd} because Theorem~\ref{thm:bernstein-subexp} allows the random matrices to be unbounded but it also demands that we control the fluctuation of the maximum \emph{and} minimum eigenvalues.  The resulting tail bound is very similar to Theorem~\ref{thm:bernstein-bdd}(ii).  We cannot achieve a Bennett-type inequality, like Theorem~\ref{thm:bernstein-bdd}(i), without stricter assumptions on the growth of moments.

The proofs of Theorem~\ref{thm:bernstein-bdd} and~\ref{thm:bernstein-subexp} appear below.  We finish the discussion with an assorted collection of enriching comments.

\begin{rem}[Rectangular Versions] \label{rem:rect-bernstein}
The matrix Bernstein inequalities admit rectangular variants.  For example, consider a sequence $\{\mtx{Z}_k\}$ of $d_1 \times d_2$ random matrices that satisfy the assumptions
$$
\Expect \mtx{Z}_k = \mtx{0}
\quad\text{and}\quad
\norm{ \mtx{Z}_k } \leq R \quad\text{almost surely}.
$$
We can apply Theorem~\ref{thm:bernstein-bdd} to the s.a.~dilation~\eqref{eqn:sa-dilation} of the sum of these random matrices to see that the probability
$$
\Prob{ \norm{ \sum\nolimits_k \mtx{Z}_k } \geq t }
	\leq d \cdot \exp\left( \frac{\sigma^2}{R^2} \cdot h\!\left( \frac{Rt}{\sigma^2} \right) \right)
$$
where $d := d_1 + d_2$ and where the variance parameter
$$
\sigma^2 := \max\left\{ \norm{\sum\nolimits_k \Expect(\mtx{Z}_k\mtx{Z}_k^\adj)}, \
	\norm{ \sum\nolimits_k \Expect( \mtx{Z}_k^\adj \mtx{Z}_k ) } \right\}.
$$
This argument leads to Theorem~\ref{thm:intro-bernstein-rect}, stated in the introduction.  There is also a rectangular extension of Theorem~\ref{thm:bernstein-subexp}, but the hypotheses are messier.
\end{rem}

\begin{rem}[Related Inequalities] 
There are too many variants of the scalar Bernstein inequality to present the matrix generalization of each one.  Let us just mention a few of the possibilities.

\begin{itemize}
\item	Theorem~\ref{thm:bernstein-subexp} can be sharpened using an idea of Rio that appears in~\cite[Sec.~2.2.3]{Mas07:Concentration-Inequalities}.

\item	When the random matrices exhibit moment growth of the form $\Expect( \mtx{X}_k^p ) \psdle R^{p-2} \mtx{A}_k^2$, we recover the Poissonian tail behavior captured in Theorem~\ref{thm:bernstein-bdd}(i).

\item	When the summands are symmetric random variables (i.e., $\mtx{X}_k \sim - \mtx{X}_k$), we can exploit the fact that the matrix mgf $\Expect \econst^{\theta \mtx{X}_k} = \Expect \cosh(\theta \mtx{X}_k)$ to obtain arcsinh inequalities.

\end{itemize}
\end{rem}


\begin{rem}[Expectations]
We can use the matrix Bernstein inequality to bound the mean of the maximum eigenvalue of the random sum.  For example, assume that the hypotheses of Theorem~\ref{thm:bernstein-bdd} or~\ref{thm:bernstein-subexp} are in force.  Then 
\begin{equation} \label{eqn:bernstein-mean}
\Expect \lambda_{\max}\left( \sum\nolimits_k \mtx{X}_k \right)
	\leq \cnst{C} \cdot \max\left\{ \sigma \sqrt{\log d}, \ R \log d \right\}.
\end{equation}
The upper bound follows by integrating Theorem~\ref{thm:bernstein-bdd}(ii) or Theorem~\ref{thm:bernstein-subexp}(i).  Lower bounds seem to require additional assumptions.
\end{rem}


\begin{rem}[Previous Work]

Oliveira's results are quite similar to the bounds presented here.  In particular, Oliveira's martingale inequality~\cite[Thm.~1.2]{Oli10:Concentration-Adjacency} implies a weaker version of Theorem~\ref{thm:bernstein-bdd}(ii).  The main result from~\cite{Oli10:Sums-Random} has a similar flavor.

\end{rem}

\subsection{Proof of Theorem~\ref{thm:bernstein-bdd}}

The main lemma shows how to bound the mgf of a zero-mean random matrix using a bound for its largest eigenvalue.

\begin{lemma}[Bounded Bernstein mgf] \label{lem:bernstein-bdd-mgf}
Suppose that $\mtx{X}$ is a random s.a.~matrix that satisfies
$$
\Expect \mtx{X} = \mtx{0}
\quad\text{and}\quad
\lambda_{\max}(\mtx{X}) \leq 1.
$$
Then
$$
\Expect \econst^{\theta \mtx{X}}
	\psdle \exp\left( (\econst^{\theta} - \theta - 1)
	\cdot \Expect(\mtx{X}^2) \right)
\quad\text{for $\theta > 0$.}
$$
\end{lemma}

As usual, the proof of the mgf bound parallels a classical method, which we learned from correspondence with Yao-Liang Yu.

\begin{proof}
Fix the parameter $\theta > 0$, and define a smooth function $f$ on the real line:
$$
f(x) = \frac{\econst^{\theta x} - \theta x - 1}{x^2}
\quad\text{for $x \neq 0$}\quad\text{and}\quad
f(0) = \frac{\theta^2}{2}.
$$
An exercise in differential calculus verifies that $f$ is increasing.  Therefore, $f(x) \leq f(1)$ when $x \leq 1$.  The eigenvalues of $\mtx{X}$ do not exceed one, so the transfer rule~\eqref{eqn:scalar-matrix} implies that
$$
f(\mtx{X}) \psdle f(1) \cdot \Id.
$$
Expanding the matrix exponential and applying the latter relation, we discover that
$$
\econst^{\theta \mtx{X}}
	= \Id + \theta \mtx{X} + \mtx{X} \cdot f(\mtx{X}) \cdot \mtx{X}
	\psdle \Id + \theta \mtx{X} + f(1) \cdot \mtx{X}^2.
$$
To complete the proof, we take the expectation of this semidefinite bound.
$$
\Expect \econst^{\theta \mtx{X}}
	\psdle \Id + f(1) \cdot \Expect( \mtx{X}^2 )
	\psdle \exp\left( f(1) \cdot \Expect( \mtx{X}^2 ) \right)
	= \exp\left( (\econst^\theta - \theta - 1) \cdot \Expect(\mtx{X}^2) \right).
$$
The second semidefinite relation follows from~\eqref{eqn:Id+A}.  
\end{proof}

We are prepared to establish the Bernstein inequalities for bounded random matrices.


\begin{proof}[Proof of Theorem~\ref{thm:bernstein-bdd}]
We assume that $R = 1$; the general result follows by a scaling argument once we note that the summands are 1-homogeneous and the variance $\sigma^2$ is 2-homogeneous.

The main challenge is to establish the Bennett inequality, Part (i); the remaining bounds are consequences of simple numerical estimates.
Invoke Lemma~\ref{lem:bernstein-bdd-mgf} to see that
$$
\Expect \econst^{\theta \mtx{X}_k}
	\psdle \exp\left( g(\theta) \cdot \Expect \big(\mtx{X}_k^2 \big) \right)
\quad\text{where $g(\theta) := \econst^{\theta} - \theta - 1$ for $\theta > 0$}.
$$
For each $\theta > 0$, Corollary~\ref{cor:main-result-indep} implies that
\begin{align*}
\Prob{ \lambda_{\max}\left( \sum\nolimits_k \mtx{X}_k \right) \geq t }
	&\leq d \cdot \exp\left( -\theta t + g(\theta) \cdot \lambda_{\max}\left(
	\sum\nolimits_k \Expect \big(\mtx{X}_k^2 \big) \right) \right) \\
	&= d \cdot \exp\left( -\theta t + g(\theta) \cdot \sigma^2 \right).
\end{align*}
The right-hand side attains its minimal value when $\theta = \log(1 + t / \sigma^2)$.  Substitute and simplify to establish Part (i).

The Bennett inequality (i) implies the Bernstein inequality (ii) because of the numerical bound
$$
h(u) \geq \frac{u^2/2}{1 + u/3}
\quad\text{for $u \geq 0$.}
$$
The latter relation is established by comparing derivatives.

The Bernstein inequality (ii) implies the split Bernstein inequality (iii).  To obtain the subgaussian piece of (iii), observe that
$$
\frac{1}{\sigma^2 + Rt/3} \geq \frac{1}{\sigma^2 + R (\sigma^2/R) / 3}
= \frac{3}{4\sigma^2}
\quad\text{for $t \leq \sigma^2/R$}
$$
because the left-hand side is a decreasing function of $t$ for $t \geq 0$.  Similarly, we obtain the subexponential piece of (iii) from the fact
$$
\frac{t}{\sigma^2 + Rt/3} \geq \frac{(\sigma^2/R)}{\sigma^2 + R (\sigma^2/R) / 3}
= \frac{3}{4R}
\quad\text{for $t \geq \sigma^2/R$},
$$
which holds because the left-hand side is an increasing function of $t$ for $t \geq 0$.
\end{proof} 

\subsection{Proof of Theorem~\ref{thm:bernstein-subexp}}

We begin with the appropriate estimate for the matrix mgf.

\begin{lemma}[Subexponential Bernstein mgf] \label{lem:bernstein-subexp-mgf}
Suppose that $\mtx{X}$ is a random s.a.~matrix that satisfies
$$
\Expect \mtx{X} = \mtx{0}
\quad\text{and}\quad
\Expect( \mtx{X}^p ) \psdle \frac{p!}{2} \cdot \mtx{A}^2
\quad\text{for $p = 2, 3, 4, \dots$.}
$$
Then
$$
\Expect \econst^{\theta \mtx{X}}
	\psdle \exp\left( \frac{\theta^2}{2(1-\theta)} \cdot \mtx{A}^2 \right)
\quad\text{for $0 < \theta < 1$.}
$$
\end{lemma}

\begin{proof}
The argument proceeds by estimating each term in the Taylor series of the matrix exponential.  Indeed,
\begin{multline*}
\Expect \econst^{\theta \mtx{X}}
	= \Id + \theta \Expect \mtx{X} + \sum_{p=2}^{\infty} \frac{\theta^{p} \Expect( \mtx{X}^{p} )}{p!}
	\psdle \Id + \sum_{p=2}^{\infty} \frac{\theta^p}{2} \cdot \mtx{A}^2
	= \Id + \frac{\theta^2}{2(1 - \theta)} \cdot \mtx{A}^2
	\psdle \exp\left( \frac{\theta^2}{2(1 - \theta)} \cdot \mtx{A}^2 \right).
\end{multline*}
As usual, the last relation is~\eqref{eqn:Id+A}.
\end{proof}

The Bernstein inequality for subexponential random matrices is an easy consequence of the previous lemma.

\begin{proof}[Proof of Theorem~\ref{thm:bernstein-subexp}]
As before, we assume that $R = 1$; the general result follows by scaling.  Invoke Lemma~\ref{lem:bernstein-subexp-mgf} to see that
$$
\Expect \econst^{\theta \mtx{X}_k}
	\psdle \exp\left( g(\theta) \cdot \mtx{A}_k^2 \right)
\quad\text{where}\quad
g(\theta) := \frac{\theta^2}{2(1 - \theta)}
\quad\text{for $0 < \theta < 1$}.
$$
For each $\theta > 0$, Corollary~\ref{cor:main-result-indep} implies that
$$
\Prob{ \lambda_{\max}\left( \sum\nolimits_k \mtx{X}_k \right) \geq t }
	\leq d \cdot \exp\left( -\theta t + g(\theta) \cdot \lambda_{\max}\left(
	\sum\nolimits_k \mtx{A}_k^2 \right) \right)
	= d \cdot \exp\left( -\theta t + g(\theta) \cdot \sigma^2 \right).
$$
We select $\theta = t/(\sigma^2 + t)$.  Substitute and simplify to complete Part (i).

The split inequality (ii) follows from Part (i) by the same argument presented in the proof of Theorem~\ref{thm:bernstein-bdd}.
\end{proof}

\section{The Matrix Hoeffding, Azuma, and McDiarmid Inequalities} \label{sec:azuma}

In this section, we prove some simple martingale deviation bounds by modifying the approach that we have used to study sums of independent random matrices.  More sophisticated martingale results require additional machinery~\cite{Oli10:Concentration-Adjacency,Tro11:Freedmans-Inequality}.

\subsection{Matrix Martingales}

We begin with the required definitions.
Let $(\Omega, \coll{F}, \mathbb{P})$ be a master probability space.  Consider a filtration $\{\coll{F}_k\}$ contained in the master sigma algebra:
$$
\coll{F}_0 \subset \coll{F}_1 \subset \coll{F}_2 \subset \dots
	\subset \coll{F}_\infty \subset \coll{F}.
$$
Given such a filtration, we define the conditional expectation
$
\Expect_k[\;\cdot\;] := \Expect[ \; \cdot \; | \; \coll{F}_k ].
$
A sequence $\{ \mtx{X}_k \}$ of random matrices is \emph{adapted} to the filtration when each $\mtx{X}_k$ is measurable with respect to $\coll{F}_k$.  Loosely speaking, an adapted sequence is one where the present depends only upon the past.

An adapted sequence $\{ \mtx{Y}_k \}$ of s.a.~matrices is called a \term{matrix martingale} when
$$
\Expect_{k-1} \mtx{Y}_k  = \mtx{Y}_{k-1}
\quad\text{and}\quad
\Expect \norm{ \mtx{Y}_k } < \infty
\quad\text{for $k = 1, 2, 3, \dots$.}
$$
We obtain a scalar martingale if we track any fixed coordinate of a matrix martingale $\{\mtx{Y}_k\}$.
Given a matrix martingale $\{\mtx{Y}_k\}$, we can construct the \term{difference sequence}
$$
\mtx{X}_k := \mtx{Y}_k - \mtx{Y}_{k-1}
\quad\text{for $k = 1, 2, 3, \dots$.}
$$
Note that the difference sequence is conditionally zero mean: $\Expect_{k-1} \mtx{X}_k  = \mtx{0}$.

\subsection{Main Results}

The scalar version of Azuma's inequality states that a scalar martingale exhibits normal concentration about its mean value, and the scale for deviations is controlled by the total maximum squared range of the difference sequence.  Here is a matrix extension.


\begin{thm}[Matrix Azuma] \label{thm:matrix-azuma}
Consider a finite adapted sequence $\{ \mtx{X}_k \}$ of self-adjoint matrices in dimension $d$, and a fixed sequence $\{ \mtx{A}_k \}$ of self-adjoint matrices that satisfy
$$
\Expect_{k-1} \mtx{X}_k = \mtx{0}
\quad\text{and}\quad
\mtx{X}_k^2 \psdle \mtx{A}_k^2
\quad\text{almost surely}.
$$
Compute the variance parameter
\begin{equation} \label{eqn:maxsq-range}
\sigma^2 := \norm{ \sum\nolimits_k \mtx{A}_k^2 }.
\end{equation}
Then, for all $t \geq 0$,
\begin{equation} \label{eqn:matrix-azuma}
\Prob{ \lambda_{\max}\left( \sum\nolimits_k \mtx{X}_k \right) \geq t }
	\leq d \cdot \econst^{-t^2/8\sigma^2}.
\end{equation}
\end{thm}

Theorem~\ref{thm:matrix-azuma} can also be phrased directly in terms of a matrix martingale.

\begin{cor}
Consider an s.a.~matrix martingale $\{\mtx{Y}_k : k = 1, \dots, n \}$ in dimension $d$, and let $\{\mtx{X}_k\}$ be the associated difference sequence.  Suppose that the difference sequence satisfies the hypotheses of Theorem~\ref{thm:matrix-azuma}, and compute the parameter $\sigma^2$ according to~\eqref{eqn:maxsq-range}.  Then
\begin{equation} \label{eqn:azuma-martingale}
\Prob{ \lambda_{\max}(\mtx{Y}_n - \Expect \mtx{Y}_n) \geq t }
	\leq d \cdot \econst^{-t^2/8\sigma^2}.
\end{equation}
\end{cor}

We continue with a few tangential comments.

\begin{rem}[Rectangular Version]
The matrix Azuma inequality has a rectangular version, which we obtain by applying Theorem~\ref{thm:matrix-azuma} to the s.a.~dilation~\eqref{eqn:sa-dilation} of the adapted sequence.
\end{rem}

\begin{rem}[Related Inequalities]
There are several situations where the constant 1/8 in the bound~\eqref{eqn:matrix-azuma} can be improved to 1/2.  One case occurs when each summand $\mtx{X}_k$ is conditionally symmetric; see Remark~\ref{rem:cond-sym}.  Another example requires the assumption that $\mtx{X}_k$ commutes almost surely with $\mtx{A}_k$, which allows us to generalize the classical proof~\cite[Lem.~2.6]{McD98:Concentration} of the Azuma inequality to the matrix setting.

If we place the additional assumption that the summands are independent, Theorem~\ref{thm:matrix-azuma} gives a matrix extension of one of Hoeffding's inequalities, which we have presented as Theorem~\ref{thm:intro-hoeffding} in the introduction.  
\end{rem}

In the scalar setting, one of the most useful corollaries of Azuma's inequality is the bounded differences inequality of McDiarmid~\cite[Thm.~3.1]{McD98:Concentration}.  This result states that a function of independent random variables exhibits normal concentration about its mean, and the variance depends on how much a change in a single variable can alter the value of the function.  A version of the bounded differences inequality holds in the matrix setting.

\begin{cor}[Matrix Bounded Differences] \label{cor:mcdiarmid}
Let $\{ Z_k : k = 1, 2, \dots, n \}$ be an independent family of random variables, and let $\mtx{H}$ be a function that maps $n$ variables to a self-adjoint matrix of dimension $d$.  Consider a sequence $\{\mtx{A}_k\}$ of fixed self-adjoint matrices that satisfy
$$
\left( \mtx{H}(z_1, \dots, z_k, \dots, z_n) - \mtx{H}(z_1, \dots, z_k', \dots, z_n ) \right)^2 \psdle \mtx{A}_k^2,
$$
where $z_i$ and $z_i'$ range over all possible values of $Z_i$ for each index $i$.  Compute the variance parameter
$$
\sigma^2 := \norm{ \sum\nolimits_k \mtx{A}_k^2 }.
$$
Then, for all $t \geq 0$,
$$
\Prob{ \lambda_{\max}( \mtx{H}(\vct{z}) - \Expect \mtx{H}(\vct{z}) ) \geq t }
	\leq d \cdot \econst^{-t^2/8\sigma^2}
$$
where $\vct{z} = (Z_1, \dots, Z_n)$.
\end{cor}

The proofs of the matrix Azuma and McDiarmid inequalities appear in the next two sections.

\subsection{Proof of Theorem~\ref{thm:matrix-azuma}}

The classical approach to Azuma's inequality does not seem to extend directly to the matrix setting.  See~\cite[Lem.~2.6]{McD98:Concentration} for a short presentation of this argument.  We use a different type of proof that is inspired by methods from probability in Banach space~\cite{LT91:Probability-Banach}.  The main idea is to inject additional randomness into the sum via a symmetrization procedure.




\begin{lemma}[Symmetrization] \label{lem:symmetrization}
Let $\mtx{H}$ be a fixed s.a.~matrix, and let $\mtx{X}$ be a random s.a.~matrix with $\Expect \mtx{X} = \mtx{0}$.  Then
$$
\Expect \trace \econst^{\mtx{H} + \mtx{X}}
	\leq \Expect \trace \econst^{\mtx{H} + 2\eps\mtx{X}},
$$
where $\eps$ is a Rademacher variable independent from $\mtx{X}$.
\end{lemma}

\begin{proof}
Construct an independent copy $\mtx{X}'$ of the random matrix, and let $\Expect'$ denote integration with respect to the new variable.  Since the matrix is zero mean,
$$
\Expect \trace \econst^{\mtx{H} + \mtx{X}}
	= \Expect \trace \econst^{\mtx{H} + \mtx{X} - \Expect' \mtx{X}'}
	\leq \Expect \trace \econst^{\mtx{H} + (\mtx{X} - \mtx{X}')}
	= \Expect \trace \econst^{\mtx{H} + \eps(\mtx{X} - \mtx{X}')}.
$$
We have used the convexity of the trace exponential to justify Jensen's inequality.
Since $\mtx{X} - \mtx{X}'$ is a symmetric random variable, we can modulate it by an independent Rademacher variable $\eps$ without changing its distribution.  The final bound depends on a short sequence of inequalities:
\begin{multline*}
\Expect \trace \econst^{\mtx{H} + \mtx{X}}
	\leq \Expect \trace \left(
	\econst^{ \mtx{H}/2 + \eps \mtx{X} } \cdot
	\econst^{ \mtx{H}/2 - \eps \mtx{X}' } \right)
	\leq \Expect \left[ \left( \trace \econst^{\mtx{H} + 2\eps\mtx{X}} \right)^{1/2} \cdot
	\big( \trace \econst^{\mtx{H} - 2\eps\mtx{X}'} \big)^{1/2} \right] \\
	\leq \left(\Expect \trace \econst^{\mtx{H} + 2\eps\mtx{X}} \right)^{1/2} \cdot
	\big( \Expect \trace \econst^{\mtx{H} - 2\eps\mtx{X}'} \big)^{1/2}
	= \Expect \trace \econst^{\mtx{H} + 2\eps\mtx{X}}.
\end{multline*}
The first relation is the Golden--Thompson inequality~\eqref{eqn:golden-thompson}; the second is the Cauchy--Schwarz inequality 
for the trace; and the third  is the Cauchy--Schwarz inequality 
for real random variables.  The last identity follows because the two factors are identically distributed.
\end{proof}

The other essential ingredient in the proof is a conditional bound for the matrix cgf of a symmetrized random matrix.

\begin{lemma}[Azuma cgf] \label{lem:azuma-cgf}
Suppose that $\mtx{X}$ is a random s.a.~matrix and $\mtx{A}$ is a fixed s.a.~matrix that satisfy $\mtx{X}^2 \psdle \mtx{A}^2$.  Let $\eps$ be a Rademacher random variable independent from $\mtx{X}$.  Then
$$
\log \Expect\big[ \econst^{2\eps \theta \mtx{X}} \, | \, \mtx{X} \big]
	\psdle 2\theta^2 \mtx{A}^2
	\quad\text{for $\theta \in \mathbb{R}$.}
$$
\end{lemma}

\begin{proof}
We apply the Rademacher mgf bound, Lemma~\ref{lem:rad-gauss-mgf}, conditionally to obtain
$$
\Expect\big[ \econst^{2\theta \eps \mtx{X}} \, | \, \mtx{X} \big]
	\psdle \econst^{2\theta^2 \mtx{X}^2}.
$$
The fact~\eqref{eqn:log-monotone} that the logarithm is operator monotone implies that
$$
\log \Expect\big[ \econst^{2\theta \eps \mtx{X}} \, | \, \mtx{X} \big]
	\psdle 2\theta^2 \mtx{X}^2
	\psdle 2\theta^2 \mtx{A}^2,
$$
where the second relation follows from the hypothesis on $\mtx{X}$.
\end{proof}

We are prepared to establish the matrix Azuma inequality.  The proof involves an iteration similar to the argument that implies the subadditivity of cgfs, Lemma~\ref{lem:cgf-indep}, for sums of independent random matrices.

\begin{proof}[Proof of Theorem~\ref{thm:matrix-azuma}]
The matrix Laplace transform method, Proposition~\ref{prop:laplace-transform}, states that
\begin{equation} \label{eqn:azuma-lt}
\Prob{ \lambda_{\max}\left( \sum\nolimits_k \mtx{X}_k \right) \geq t }
	\leq \inf_{\theta > 0} \left\{ \econst^{-\theta t} \cdot
	\Expect \trace \exp\left( \sum\nolimits_k \theta \mtx{X}_k \right) \right\}.
\end{equation}
The main difficulty in the proof is to bound the matrix mgf, which we accomplish by an iterative argument that alternates between symmetrization and cumulant bounds.

Let us detail the first step of the iteration.  Define the natural filtration $\coll{F}_{k} := \coll{F}(\mtx{X}_1, \dots, \mtx{X}_k)$ of the process $\{\mtx{X}_k\}$.  Then we may compute
\begin{align*}
\Expect \trace \exp\left( \sum\nolimits_k \theta \mtx{X}_k \right)
	&= \Expect \Expect\left[ \trace \exp\left( \sum\nolimits_{k=1}^{n-1}
		\theta \mtx{X}_k + \theta \mtx{X}_n \right)
		\, \big| \, \coll{F}_{n-1} \right] \\
	&\leq \Expect \Expect\left[ \trace \exp\left(
	\sum\nolimits_{k=1}^{n-1} \theta\mtx{X}_k
	+ 2\eps \theta \mtx{X}_n \right) \, \big| \, \coll{F}_{n} \right] \\
	&\leq \Expect \trace \exp\left( \sum\nolimits_{k=1}^{n-1} \theta \mtx{X}_k
	+ \log \Expect\big[ \econst^{2\eps \theta \mtx{X}_n} \, | \, \coll{F}_n \big] \right) \\
	&\leq \Expect \trace \exp\left( \sum\nolimits_{k=1}^{n-1} \theta \mtx{X}_k
	+ 2\theta^2 \mtx{A}_n^2 \right).
\end{align*}
The first identity is the tower property of conditional expectation.  In the second line, we invoke the symmetrization method, Lemma~\ref{lem:symmetrization}, conditional on $\coll{F}_{n-1}$, and then we relax the conditioning on the inner expectation to the larger algebra $\coll{F}_n$.  By construction, the Rademacher variable $\eps$ is independent from $\coll{F}_n$, so we can apply the concavity result, Corollary~\ref{cor:cum-ineq}, conditional on $\coll{F}_n$.  Finally, we use the fact~\eqref{eqn:exp-trace-monotone} that the trace exponential is monotone to introduce the Azuma cgf bound, Lemma~\ref{lem:azuma-cgf}, in the last inequality.

By iteration, we achieve
\begin{equation} \label{eqn:azuma-mgf-pf}
\Expect \trace \exp\left( \sum\nolimits_k \theta \mtx{X}_k \right)
	\leq \trace \exp\left( 2\theta^2 \sum\nolimits_k \mtx{A}_k^2 \right).
\end{equation}
Note that this procedure relies on the fact that the sequence $\{\mtx{A}_k\}$ of upper bounds does not depend on the values of the random sequence $\{\mtx{X}_k\}$.  Substitute the mgf bound~\eqref{eqn:azuma-mgf-pf} into the Laplace transform bound~\eqref{eqn:azuma-lt}, and observe that the infimum is achieved when $\theta = t/4\sigma^2$.
\end{proof}

\begin{rem} \label{rem:cond-sym}
Suppose that the sequence $\{\mtx{X}_k\}$ is conditionally symmetric:
$$
\mtx{X}_k \sim - \mtx{X}_k
\quad\text{conditional on $\coll{F}_{k-1}$.}
$$
When we execute the proof of Theorem~\ref{thm:matrix-azuma} under this assumption, we can symmetrize each term in the sum without suffering an extra factor of two.  For example,
$$
\Expect\left[ \trace \exp\left( \sum\nolimits_{k=1}^{n-1} \theta \mtx{X}_k
	+ \theta \mtx{X}_n \right) \, \big| \, \coll{F}_{n-1} \right]
	= \Expect\left[ \trace \exp\left( \sum\nolimits_{k=1}^{n-1} \theta \mtx{X}_k
	+ \eps \theta \mtx{X}_n \right) \, \big| \, \coll{F}_{n-1} \right]
$$
where $\eps$ is independent from $\coll{F}_n$.  The rest of the proof remains the same, but the analog of the bound~\eqref{eqn:matrix-azuma} has a constant of 1/2 instead of 1/8 in the exponent. 
\end{rem}

\subsection{Proof of Corollary~\ref{cor:mcdiarmid}}

Finally, we establish the matrix version of the bounded differences inequality.  The main idea in the argument is to construct the Doob martingale associated with the natural filtration of the independent random sequence.  We compute semidefinite bounds for the difference sequence, and then we apply the matrix Azuma inequality to control the deviations of the martingale.

\begin{proof}[Proof of Corollary~\ref{cor:mcdiarmid}]
In this argument only, we write $\Expect_{Z}$ for the expectation with respect to a random variable $Z$, holding other variables fixed.  Recall that $\vct{z} = (Z_1, \dots, Z_n)$.  For $k = 0, 1, \dots, n$, consider the random matrices
$$
\mtx{Y}_k := \Expect[ \mtx{H}(\vct{z}) \; | \; Z_1, Z_2, \dots, Z_{k} ]
	= \Expect_{Z_{k+1}} \Expect_{Z_{k+2}} \dots \Expect_{Z_n}  \mtx{H}(\vct{z}) .
$$
The sequence $\{\mtx{Y}_k\}$ forms a Doob martingale.  The associated difference sequence is
$$
\mtx{X}_k
	:= \mtx{Y}_k - \mtx{Y}_{k-1}
	= \Expect_{Z_{k+1}} \Expect_{Z_{k+2}} \dots \Expect_{Z_n} \big( \mtx{H}(\vct{z}) - \Expect_{Z_{k}} \mtx{H}(\vct{z}) \big),
$$
where the second identity follows from independence and Fubini's theorem.

It remains to bound the difference sequence.
Let $Z_k'$ be an independent copy of $Z_k$, and construct the random vector $\vct{z}' = (Z_1, \dots, Z_{k-1}, Z_k', Z_{k+1}, \dots, Z_n)$.  Observe that
$
\Expect_{Z_k} \mtx{H}(\vct{z}) = \Expect_{Z_k'} \mtx{H}(\vct{z}')
$
and that $\mtx{H}(\vct{z})$ does not depend on $Z_k'$.  Therefore, we can write
$$
\mtx{X}_k = \Expect_{Z_{k+1}} \Expect_{Z_{k+2}} \dots \Expect_{Z_n} \Expect_{Z_k'} \big(\mtx{H}(\vct{z})-\mtx{H}(\vct{z}') \big).
$$
The vectors $\vct{z}$ and $\vct{z}'$ differ only in the $k$th coordinate, so that
$$
\big(\mtx{H}(\vct{z}) - \mtx{H}(\vct{z}') \big)^2 \psdle \mtx{A}_k^2
$$
by definition of the bound $\mtx{A}_k^2$.  Finally, the semidefinite Jensen inequality~\eqref{eqn:jensen-matrix-square} for the matrix square yields
$$
\mtx{X}_k^2 \psdle \Expect_{Z_{k+1}} \Expect_{Z_{k+2}} \dots \Expect_{Z_n} \Expect_{Z_k'} \big( \mtx{H}(\vct{z}) - \mtx{H}(\vct{z}') \big)^2
	\psdle \mtx{A}_k^2.
$$
To complete the proof, we apply~\eqref{eqn:azuma-martingale} to the martingale $\{\mtx{Y}_k\}$.
\end{proof}

\section*{Acknowledgments}

I would like to thank Vern Paulsen and Bernhard Bodmann for some helpful conversations connected with this project.  Klas Markstr{\"o}m and David Gross provided references to related work.  Ben Recht offered some useful comments on the presentation.  Yao-Liang Yu proposed the argument in Lemma~\ref{lem:bernstein-bdd-mgf}.  Richard Chen and Alex Gittens have helped me root out (numerous) typographic errors.  Finally, let me mention Roberto Oliveira's elegant work~\cite{Oli10:Sums-Random} on matrix probability inequalities, which originally spurred me to pursue this project.  





\bibliographystyle{alpha}
\bibliography{user-friendly}

\end{document}

%% file: macro-file.tex
\usepackage{amsmath}
\usepackage{amssymb}
\usepackage{bm}
\usepackage{mathrsfs}
\usepackage[dvips]{graphicx}

\usepackage{color}

\usepackage{url}

\usepackage{supertabular}

\usepackage{alg}

\makeatletter
\def\algspacing{\alg@unmargin}
\makeatother

\newtheorem{thm}{Theorem}
\newtheorem{cor}[thm]{Corollary}
\newtheorem{lemma}[thm]{Lemma}
\newtheorem{prop}[thm]{Proposition}

\theoremstyle{remark}
\newtheorem{rem}[thm]{Remark}

\numberwithin{equation}{section}

\numberwithin{thm}{section}

\DeclareMathAlphabet{\mathsfsl}{OT1}{cmss}{m}{sl}

\newcommand{\lang}{\textit}

\newcommand{\term}{\emph}

\newcommand{\cnst}[1]{\mathrm{#1}}

	{\makebox{\phantom{}} \\ %
		\textsc{Input:}\begin{itemize}}%
	{\end{itemize}}

	{\textsc{Output:}\begin{itemize}}%
	{\end{itemize}}

	{\textsc{Procedure:}\begin{enumerate}}%
	{\end{enumerate}}

\renewcommand{\phi}{\varphi}

\newcommand{\eps}{\varepsilon}

\newcommand{\econst}{\mathrm{e}}

\newcommand{\onevct}{\mathbf{e}}

\newcommand{\Id}{\mathbf{I}}

\newcommand{\coll}[1]{\mathscr{#1}}

\newcommand{\Rspace}[1]{\mathbb{R}^{#1}}

\newcommand{\abs}[1]{\left\vert {#1} \right\vert}
\newcommand{\abssq}[1]{{\abs{#1}}^2}

\newcommand{\diff}[1]{\mathrm{d}{#1}}
\newcommand{\idiff}[1]{\, \diff{#1}}

\newcommand{\Prob}[1]{\mathbb{P}\left\{ {#1} \right\}}
\newcommand{\Expect}{\operatorname{\mathbb{E}}}

\newcommand{\vct}[1]{\bm{#1}}
\newcommand{\mtx}[1]{\bm{#1}}

\newcommand{\adj}{*}

\newcommand{\diag}{\operatorname{diag}}
\newcommand{\trace}{\operatorname{tr}}

\newcommand{\psdle}{\preccurlyeq}
\newcommand{\psdge}{\succcurlyeq}

\newcommand{\psdlt}{\prec}

\newcommand{\norm}[1]{\left\Vert {#1} \right\Vert}
\newcommand{\normsq}[1]{\norm{#1}^2}

\newcommand{\smnorm}[2]{{\bigl\Vert {#2} \bigr\Vert}_{#1}}

\newcommand{\enorm}[1]{\norm{#1}_2}
\newcommand{\enormsq}[1]{\enorm{#1}^2}






%% file: User-Friendly-focm-vf.bbl
\begin{thebibliography}{HMT11}

\bibitem[AC09]{AC09:Fast-Johnson-Lindenstrauss}
N.~Ailon and B.~Chazelle.
\newblock The fast {J}ohnson--{L}indenstrauss transform and approximate nearest
  neighbors.
\newblock {\em SIAM J. Comput.}, 39(1):302--322, 2009.

\bibitem[AM07]{AM07:Fast-Computation}
D.~Achlioptas and F.~McSherry.
\newblock Fast computation of low-rank matrix approximations.
\newblock {\em J. Assoc. Comput. Mach.}, 54(2):Article 10, 2007.
\newblock (electronic).

\bibitem[AW02]{AW02:Strong-Converse}
R.~Ahlswede and A.~Winter.
\newblock Strong converse for identification via quantum channels.
\newblock {\em IEEE Trans. Inform. Theory}, 48(3):569--579, Mar. 2002.

\bibitem[Bha97]{Bha97:Matrix-Analysis}
R.~Bhatia.
\newblock {\em Matrix Analysis}.
\newblock Number 169 in Graduate Texts in Mathematics. Springer, Berlin, 1997.

\bibitem[Bha07]{Bha07:Positive-Definite}
R.~Bhatia.
\newblock {\em Positive Definite Matrices}.
\newblock Princeton Univ. Press, Princeton, NJ, 2007.

\bibitem[Bog98]{Bog98:Gaussian-Measures}
V.~Bogdanov.
\newblock {\em Gaussian Measures}.
\newblock American Mathematical Society, Providence, RI, 1998.

\bibitem[Buc01]{Buc01:Operator-Khintchine}
A.~Buchholz.
\newblock Operator {K}hintchine inequality in non-commutative probability.
\newblock {\em Math. Ann.}, 319:1--16, 2001.

\bibitem[Buc05]{Buc05:Optimal-Constants}
A.~Buchholz.
\newblock Optimal constants in {K}hintchine-type inequalities for {F}ermions,
  {R}ademachers and $q$-{G}aussian operators.
\newblock {\em Bull. Pol. Acad. Sci. Math.}, 53(3):315--321, 2005.

\bibitem[Che52]{Che52:Measure-Asymptotic}
H.~Chernoff.
\newblock A measure of the asymptotic efficiency for tests of a hypothesis
  based on the sum of observations.
\newblock {\em Ann. Math. Statist.}, 23(4):493--507, 1952.

\bibitem[CM08]{CM08:Expansion-Properties}
D.~Cristofides and K.~Markstr{\"o}m.
\newblock Expansion properties of random {C}ayley graphs and vertex transitive
  graphs via matrix martingales.
\newblock {\em Random Structures Algs.}, 32(8):88--100, 2008.

\bibitem[CR07]{CR07:Sparsity-Incoherence}
E.~Cand{\`e}s and J.~K. Romberg.
\newblock Sparsity and incoherence in compressive sampling.
\newblock {\em Inverse Problems}, 23(3):969--985, 2007.

\bibitem[dlPG02]{PG02:Decoupling}
V.~H. de~la Pe{\~n}a and E.~Gin{\'e}.
\newblock {\em Decoupling: From Dependence to Independence}.
\newblock Probability and its Applications. Springer, Berlin, 2002.

\bibitem[DS02]{DS02:Local-Operator}
K.~R. Davidson and S.~J. Szarek.
\newblock Local operator theory, random matrices, and {B}anach spaces.
\newblock In W.~B. Johnson and J.~Lindenstrauss, editors, {\em Handbook of
  Banach Space Geometry}, pages 317--366. Elsevier, Amsterdam, 2002.

\bibitem[DZ98]{DZ98:Large-Deviations}
A.~Dembo and O.~Zeitouni.
\newblock {\em Large Deviations: Techniques and Applications}.
\newblock Springer, 2nd edition, 1998.

\bibitem[Eff09]{Eff09:Matrix-Convexity}
E.~G. Effros.
\newblock A matrix convexity approach to some celebrated quantum inequalities.
\newblock {\em Proc. Natl. Acad. Sci. USA}, 106(4):1006--1008, Jan. 2009.

\bibitem[Eps73]{Eps73:Remarks-Two}
H.~Epstein.
\newblock Remarks on two theorems of {E.} {L}ieb.
\newblock {\em Comm. Math. Phys.}, 31:317--325, 1973.

\bibitem[Fre75]{Fre75:Tail-Probabilities}
D.~A. Freedman.
\newblock On tail probabilities for martingales.
\newblock {\em Ann. Probab.}, 3(1):100--118, Feb. 1975.

\bibitem[Gor85]{Gor85:Some-Inequalities}
Y.~Gordon.
\newblock Some inequalities for {G}aussian processes and applications.
\newblock {\em Israel J. Math.}, 50(4):265--289, 1985.

\bibitem[Gor92]{Gor92:Majorization-Gaussian}
Y.~Gordon.
\newblock Majorization of {G}aussian processes and geometric applications.
\newblock {\em Probab. Theory Related Fields}, 91(2):251--267, 1992.

\bibitem[Gro11]{Gro11:Recovering-Low-Rank}
D.~Gross.
\newblock Recovering low-rank matrices from few coefficients in any basis.
\newblock {\em IEEE Trans. Inform. Theory}, 57(3):1548--1566, Mar. 2011.

\bibitem[Hig08]{Hig08:Functions-Matrices}
N.~J. Higham.
\newblock {\em Functions of Matrices: Theory and Computation}.
\newblock Society for Industrial and Applied Mathematics, Philadelphia, PA,
  2008.

\bibitem[HJ85]{HJ85:Matrix-Analysis}
R.~A. Horn and C.~R. Johnson.
\newblock {\em Matrix Analysis}.
\newblock Cambridge Univ. Press, Cambridge, 1985.

\bibitem[HJ94]{HJ94:Topics-Matrix}
R.~A. Horn and C.~R. Johnson.
\newblock {\em Topics in Matrix Analysis}.
\newblock Cambridge Univ. Press, Cambridge, 1994.

\bibitem[HMT11]{HMT11:Finding-Structure}
N.~Halko, P.-G. Martinsson, and J.~A. Tropp.
\newblock Finding structure with randomness: {S}tochastic algorithms for
  constructing approximate matrix decompositions.
\newblock {\em SIAM Rev.}, 53(2):217--288, June 2011.

\bibitem[HP03]{HP03:Jensens-Operator}
F.~Hansen and G.~K. Pedersen.
\newblock {J}ensen's operator inequality.
\newblock {\em Bull. London Math. Soc.}, 35:553--564, 2003.

\bibitem[JX05]{JX05:Best-Constants}
M.~Junge and Q.~Xu.
\newblock On the best constants in some non-commutative martingale
  inequalities.
\newblock {\em Bull. London Math. Soc.}, 37:243--253, 2005.

\bibitem[JX08]{JX08:Noncommutative-Burkholder-II}
M.~Junge and Q.~Xu.
\newblock Noncommutative {B}urkholder/{R}osenthal inequalities {II}:
  {A}pplications.
\newblock {\em Israel J. Math.}, 167:227--282, 2008.

\bibitem[Lat05]{Lat05:Some-Estimates}
R.~Lata{\l}a.
\newblock Some estimates of norms of random matrices.
\newblock {\em Proc. Amer. Math. Soc.}, 133(5):1273--1282, 2005.

\bibitem[Lie73]{Lie73:Convex-Trace}
E.~H. Lieb.
\newblock Convex trace functions and the {W}igner--{Y}anase--{D}yson
  conjecture.
\newblock {\em Adv. Math.}, 11:267--288, 1973.

\bibitem[Lin74]{Lin74:Expectations-Entropy}
G.~Lindblad.
\newblock Expectations and entropy inequalities for finite quantum systems.
\newblock {\em Comm. Math. Phys.}, 39:111--119, 1974.

\bibitem[LP86]{L-P86:Inegalites-Khintchine}
F.~Lust-Piquard.
\newblock In{\'e}galit{\'e}s de {K}hintchine dans {$C_p$ $(1 < p < \infty)$}.
\newblock {\em C. R. Math. Acad. Sci. Paris}, 303(7):289--292, 1986.

\bibitem[LPP91]{LPP91:Noncommutative-Khintchine}
F.~Lust-Piquard and G.~Pisier.
\newblock Noncommutative {K}hintchine and {P}aley inequalities.
\newblock {\em Ark. Mat.}, 29(2):241--260, 1991.

\bibitem[LT91]{LT91:Probability-Banach}
M.~Ledoux and M.~Talagrand.
\newblock {\em Probability in Banach Spaces: Isoperimetry and Processes}.
\newblock Springer, Berlin, 1991.

\bibitem[Lug09]{Lug09:Concentration-Measure}
G.~Lugosi.
\newblock Concentration-of-measure inequalities.
\newblock Available at \url{http://www.econ.upf.edu/~lugosi/anu.pdf}, 2009.

\bibitem[Mas07]{Mas07:Concentration-Inequalities}
P.~Massart.
\newblock {\em Concentration Inequalities and Model Selection: Ecole d'Et{\'e}
  de Probabilit{\'e}s de Saint-Flour XXXIII---2003}.
\newblock Number 1896 in Lecture Notes in Mathematics (LNM). Springer, 2007.

\bibitem[McD98]{McD98:Concentration}
C.~McDiarmid.
\newblock Concentration.
\newblock In {\em Probabilistic Methods for Algorithmic Discrete Mathematics},
  number~16 in Algorithms and Combinatorics, pages 195--248. Springer, Berlin,
  1998.

\bibitem[MR95]{MR95:Randomized-Algorithms}
R.~Motwani and P.~Raghavan.
\newblock {\em Randomized Algorithms}.
\newblock Cambridge Univ. Press, Cambridge, 1995.

\bibitem[Nem07]{Nem07:Sums-Random}
A.~Nemirovski.
\newblock Sums of random symmetric matrices and quadratic optimization under
  orthogonality constraints.
\newblock {\em Math. Prog. Ser. B}, 109:283--317, 2007.

\bibitem[Oli10a]{Oli10:Concentration-Adjacency}
R.~I. Oliveira.
\newblock Concentration of the adjacency matrix and of the {L}aplacian in
  random graphs with independent edges.
\newblock Available at \url{arXiv:0911.0600}, Feb. 2010.

\bibitem[Oli10b]{Oli10:Sums-Random}
R.~I. Oliveira.
\newblock Sums of random {H}ermitian matrices and an inequality by {R}udelson.
\newblock {\em Electron. Commun. Probab.}, 15:203--212, 2010.

\bibitem[Par87]{Par87:Symmetric-Eigenvalue}
B.~N. Parlett.
\newblock {\em The Symmetric Eigenvalue Problem}.
\newblock Number~20 in Classics in Applied Mathematics. Society for Industrial
  and Applied Mathematics, Philadelphia, PA, 1987.

\bibitem[Pau02]{Pau02:Completely-Bounded}
V.~I. Paulsen.
\newblock {\em Completely Bounded Maps and Operator Algebras}.
\newblock Number~78 in Cambridge Studies in Advanced Mathematics. Cambridge
  Univ. Press, Cambridge, 2002.

\bibitem[Pet94]{Pet94:Survey-Certain}
D.~Petz.
\newblock A survey of certain trace inequalities.
\newblock In {\em Functional analysis and operator theory}, volume~30 of {\em
  Banach Center Publications}, pages 287--298, Warsaw, 1994. Polish Acad. Sci.

\bibitem[Pis03]{Pis03:Introduction-Operator}
G.~Pisier.
\newblock {\em Introduction to Operator Spaces}.
\newblock Cambridge Univ. Press, Cambridge, 2003.

\bibitem[Rec09]{Rec09:Simpler-Approach}
B.~Recht.
\newblock Simpler approach to matrix completion.
\newblock {\em J. Mach. Learn. Res.}, Oct. 2009.
\newblock To appear. Available at
  \url{http://pages.cs.wisc.edu/~brecht/papers/09.Recht.ImprovedMC.pdf}.

\bibitem[Rud99]{Rud99:Random-Vectors}
M.~Rudelson.
\newblock Random vectors in the isotropic position.
\newblock {\em J. Funct. Anal.}, 164:60--72, 1999.

\bibitem[Rus02]{Rus02:Inequalities-Quantum}
M.~B. Ruskai.
\newblock Inequalities for quantum entropy: {A} review with conditions for
  equality.
\newblock {\em J. Math. Phys.}, 43(9):4358--4375, Sep. 2002.

\bibitem[Rus05]{Rus05:Erratum-Inequalities}
M.~B. Ruskai.
\newblock Erratum: {I}nequalities for quantum entropy: {A} review with
  conditions for equality [\emph{{J}. {M}ath. {P}hys.} 43, 4358 (2002)].
\newblock {\em J. Math. Phys.}, 46(1):0199101, 2005.

\bibitem[RV07]{RV07:Sampling-Large}
M.~Rudelson and R.~Vershynin.
\newblock Sampling from large matrices: {A}n approach through geometric
  functional analysis.
\newblock {\em J. Assoc. Comput. Mach.}, 54(4):Article 21, 19 pp., Jul. 2007.
\newblock (electronic).

\bibitem[Seg00]{Seg00:Expected-Norm}
Y.~Seginer.
\newblock The expected norm of random matrices.
\newblock {\em Combin. Probab. Comput.}, 9:149--166, 2000.

\bibitem[So09]{So09:Moment-Inequalities}
A.~M.-C. So.
\newblock Moment inequalities for sums of random matrices and their
  applications in optimization.
\newblock {\em Math. Prog. Ser. A}, Dec. 2009.
\newblock (electronic).

\bibitem[SST06]{SST06:Smoothed-Analysis}
A.~Sankar, D.~A. Spielman, and S.-H. Teng.
\newblock Smoothed analysis of the condition numbers and growth factors of
  matrices.
\newblock {\em SIAM J. Matrix Anal. Appl.}, 28(2):446--476, 2006.

\bibitem[TJ74]{TJ74:Moduli-Smoothness}
N.~Tomczak-Jaegermann.
\newblock The moduli of smoothness and convexity and the {R}ademacher averages
  of trace classes {$S_p$ $(1 \leq p < \infty)$}.
\newblock {\em Studia Math.}, 50:163--182, 1974.

\bibitem[Tro08]{Tro08:Conditioning-Random}
J.~A. Tropp.
\newblock On the conditioning of random subdictionaries.
\newblock {\em Appl. Comput. Harmon. Anal.}, 25:1--24, 2008.

\bibitem[Tro10]{Tro10:Improved-Analysis}
J.~A. Tropp.
\newblock Improved analysis of the subsampled randomized {H}adamard transform.
\newblock {\em Adv. Adapt. Data Anal.}, 2010.
\newblock To appear. Available at \url{arXiv:1011.1595}.

\bibitem[Tro11a]{Tro11:Freedmans-Inequality}
J.~A. Tropp.
\newblock {F}reedman's inequality for matrix martingales.
\newblock {\em Electron. Commun. Probab.}, 16:262--270, 2011.

\bibitem[Tro11b]{Tro11:Joint-Convexity}
J.~A. Tropp.
\newblock From the joint convexity of quantum relative entropy to a concavity
  theorem of {L}ieb.
\newblock {\em Proc. Amer. Math. Soc.}, 2011.
\newblock To appear. Available at \url{arXiv:1101.1070}.

\bibitem[Tro11c]{Tro11:User-Friendly-Martingale-TR}
J.~A. Tropp.
\newblock User-friendly tail bounds for matrix martingales.
\newblock ACM Report 2011-01, California Inst. Tech., Pasadena, CA, Jan. 2011.

\bibitem[Ver09]{Ver09:Note-Sums}
R.~Vershynin.
\newblock A note on sums of independent random matrices after
  {A}hlswede--{W}inter.
\newblock Available at
  \url{http://www-personal.umich.edu/~romanv/teaching/reading-group/ahlswede-w%
inter.pdf}, 2009.

\end{thebibliography}
